\documentclass[12pt,draft]{article}

\usepackage{amsmath}
\usepackage{amsthm}
\usepackage{amssymb}
\usepackage{amsfonts}
\usepackage{latexsym}               
\usepackage{amsgen}
\usepackage[mathscr]{eucal}
\usepackage{mathrsfs}
\usepackage{bm}
\usepackage{enumerate}
\usepackage{mathptmx} 
\usepackage{fancybox}
\usepackage{color}
\usepackage{graphicx}

\newtheorem{lem}{Lemma}[section]
\newtheorem{thm}[lem]{Theorem}
\newtheorem{pro}[lem]{Proposition}
\newtheorem{rem}[lem]{Remark}

\numberwithin{equation}{section}

\begin{document}


\title{Stability of Stationary Solutions to 
the Nonisentropic Euler--Poisson System
in a Perturbed Half Space}

\author{Mingjie Li${}^1$ and Masahiro Suzuki${}^2$}

\date{%
\normalsize
${}^1$%
College of Science, Minzu University of China, Beijing 100081, P. R. China
\\[7pt]
${}^2$%
Department of Computer Science and Engineering, 
Nagoya Institute of Technology,
\\
Gokiso-cho, Showa-ku, Nagoya, 466-8555, Japan
}

\maketitle

\begin{abstract}
The main concern of this paper is to mathematically investigate the formation of 
a plasma sheath near the surface of nonplanar walls.
We study the existence and asymptotic stability of stationary solutions 
for the nonisentropic Euler-Poisson equations in a domain of which boundary is drawn by a graph, by employing a space weighted energy method. Moreover, the convergence rate of the solution toward the stationary solution is obtained, provided that the initial
perturbation belongs to the weighted Sobolev space. Because the domain is the perturbed half space, we first show the time-global solvability of  the nonisentropic Euler-Poisson equations, 
then construct stationary solutions by using the time-global solutions.
\end{abstract}

\begin{description}
\item[{\it Keywords:}]
plasma;
sheath;
Bohm criterion;
initial--boundary value problem;
long-time behavior;
convergence rate.

\item[{\it 2020 Mathematics Subject Classification:}]
76X05; 
35M12; 
35M13; 
35A01; 
35B35. 

\end{description}


\newpage

\section{Introduction}\label{S1}

We consider the nonisentropic Euler-Poisson system in the perturbed half space
\[
\Omega:=\{x=(x,x')=(x_1,x_2,x_3)\in \mathbb R^3 \, | 
\, x_1>M(x')\}
\quad \text{for $M \in \cap_{k=1}^\infty H^k(\mathbb R^2)$}.
\]
The nonisentropic Euler-Poisson system is written by
\begin{subequations}\label{eq0}
   \begin{gather}
   \rho_t+\nabla \cdot(\rho \bm{u})=0,
   \label{eq1}\\
   m\bm{u}_t +  m\left( \bm{u} \cdot \nabla \right) \bm{u}
    + R \theta \nabla (\log \rho) + R \nabla \theta = \nabla \phi,
   \label{eq2}\\
   \theta_t +   \bm{u} \cdot \nabla \theta 
    + (\gamma-1) \theta \nabla \cdot \bm{u} = 0,
   \label{eq4}\\
   \Delta \phi=\rho-e^{-\phi}, \quad t>0, \ x \in \Omega,
   \label{eq3}
   \end{gather}
where unknown functions
$\rho$, $\bm{u}=(u_1,u_2,u_3)$, $\theta$, and $-\phi$
represent the density, velocity, and temperature of the positive ions and
the electrostatic potential, respectively.
Furthermore, the mass of an ion $m>0$, the gas constant $R>0$ and the heat capacity ratio $\gamma>1$ are positive constants.
The first equation is the conservation of mass, the second one
is the equation of momentum in which the pressure gradient
and electrostatic potential gradient as well as the convection effect are taken into account,
and the third equation comes from the conservation of energy.
The fourth equation is the Poisson equation,
which governs the relation between the potential and the density of charged particles.
It is obtained by assuming the Boltzmann relation in which the electron density is 
given by $\rho_e=e^{-\phi}$. 
This assumption can be justified mathematically \cite{GGPS} even for the case that the ion flow is isentropic and the domain is the whole space.

When a material is surrounded by a plasma and the plasma contacts with its surface,
there appears a non-neutral potential region (boundary layer) between the surface and plasma, 
and a nontrivial equilibrium of the densities is achieved. 
This non-neutral region is referred as to {\it a sheath}. 
For the formation of sheath, Langmuir \cite{I.L.1} observed that 
positive ions must enter the sheath region with a sufficiently large kinetic energy. 
Bohm \cite{D.B.1} proposed the \emph{Bohm criterion} \eqref{Bohm1} below in the case of planar wall.
It states that the ion velocity at the plasma edge must exceed the ion acoustic speed.
For more details of  physicality of the sheath formation, 
we refer the reader to \cite{D.B.1, F.C.1,I.L.1, LL, K.R.1, K.R.2}.
In this paper, we investigate mathematically the formation of sheath.

We prescribe the initial and boundary conditions 
 \begin{gather}
  (\rho,\bm{u},\theta)(0,x)=(\rho_0,\bm{u}_0,\theta_{0})(x), 
  \label{ini1}
  \\
\lim_{x_1\to\infty}(\rho,u_1,u_2,u_3,\theta,\phi)(t,x_1,x')=(1,u_+,0,0,\theta_{+},0),
 \label{bc1}
 \\
  \phi(t,M(x'),x')=\phi_b \quad \text{for} \
  x'\in \mathbb R^2,
  \label{bc2}
 \end{gather}
 \end{subequations}
where $u_+<0$, $\theta_{+}>0$, and $\phi_b\in \mathbb R$ are constants.
Here $\phi_b$ represents the voltage on the boundary. Furthermore, $\lim_{x_{1}\to\infty}\phi(t,x_{1},x')=0$ means that the reference point of the electrostatic potential is located at $x_{1}=\infty$.
To solve \eqref{eq3} with this condition in the classical sense, we need $\lim_{x_{1}\to\infty}\rho(t,x_{1},x')=1$.
The unit outer normal vector of the boundary 
$\partial\Omega=\{x\in \mathbb R^3 \, | \, x_1=M(x')\}$ is represented by 
\begin{equation*}
\bm{n}(x')=(n_1,n_2,n_3)(x')
:=\left(\frac{-1}{\sqrt{1+|\nabla M|^2}},
\frac{\partial_{x_2}M}{\sqrt{1+|\nabla M|^2}},
\frac{\partial_{x_3}M}{\sqrt{1+|\nabla M|^2}} \right)(x').
\end{equation*}

We construct solutions in the region, where the following three conditions hold:
  \begin{gather}
  \inf_{x\in \Omega}\rho(t,x)>0, \quad  \inf_{x\in \Omega}\theta(t,x)>0,
  \label{po1}
  \\
\inf_{x \in \partial\Omega} \left( \frac{\bm{u}(t,x)\cdot\nabla(M(x')-x_1)}{\sqrt{1+|\nabla M(x')|^2}}-\sqrt{\frac{\gamma R\theta(t,x)}{m}}\right)>0,
  \label{super1'}
  \end{gather}
by assuming the same conditions for the initial data $(\rho_0,\bm{u}_0,\theta_{0})$:
\[
 \inf_{x\in \Omega}\rho_0(x)>0, \quad  \inf_{x\in \Omega}\theta_0(x)>0, \quad 
\inf_{x \in \partial\Omega} \left( \frac{\bm{u}_0(x)\cdot\nabla(M(x')-x_1)}{\sqrt{1+|\nabla M(x')|^2}}-\sqrt{\frac{\gamma R\theta_0(x)}{m}} \right)>0. 
\]
In particular, the supersonic outflow condition \eqref{super1'} is necessary for 
the well-posedness of the initial--boundary value problem \eqref{eq0},
because it guarantees that no boundary condition is suitable for the equations \eqref{eq1}--\eqref{eq4}.
In this case, no compatibility condition is required.
For the end state of velocity $u_+$ and temperature $\theta_{+}$, 
we assume the Bohm criterion and the supersonic outflow condition:
\begin{gather}
mu_+^2>\gamma R \theta_{+} +1, \quad u_+<0,
\label{Bohm1}\\
\inf_{x \in \partial\Omega} 
\frac{-u_+}{\sqrt{1+|\nabla M(x')|^2}}-\sqrt{\frac{\gamma R\theta_+}{m}}>0.
\label{asp1'}
\end{gather}
The thing is that we need \eqref{asp1'} to establish solutions of the problem \eqref{eq0}
in a neighborhood of the end state 
$(\rho,u_1,u_2,u_3,\theta,\phi)=(1,u_+,0,0,\theta_{+},0)$,
which is a trivial solution for the case $\phi_b=0$.

Nowadays there are many mathematical studies on the sheath formation by using the Euler--Poisson system. 
We begin by reviewing those studies for the isothermal or isentropic flow with planar walls.
Ambroso--M\'ehats--Raviart made a pioneering work \cite{AMR}.
They established the unique existence of the monotone stationary solutions
over a bounded interval, provided that the Bohm criterion \eqref{Bohm1} holds.
Moreover, it is numerically checked by Ambroso \cite{A.A.} that
the solutions of the initial--boundary value problem approaches the stationary solutions established in \cite{AMR}
as the time variable becomes large.
After that, Suzuki \cite{M.S.1} derived a necessary and sufficient condition 
for the existence of the monotone stationary solutions over a half space, 
and also found out that the Bohm criterion is a sufficient condition but not a necessary condition.
Furthermore, Ohnawa--Nishibata--Suzuki \cite{NOS} showed the asymptotic stability of the stationary solutions assuming the Bohm criterion.
The similar results were also obtained for a multicomponent plasma which consists of electrons and several components of ions \cite{M.S.2}.
These results validated mathematically the Bohm criterion
and defined the fact that the sheath corresponds to the stationary solution.

We also review the studies on the quasi-neutral limit problem
as letting the Debye length in the Euler--Poisson system tend to zero.
It is another direction to understand the sheath formation.
Ambroso--M\'ehats--Raviart \cite{AMR} and Slemrod--Sternberg \cite{SS1} investigated the problem 
over a bounded interval.
Furthermore, G\'erard-Varet--Han-Kwan--Rousset \cite{GHR1,GHR2}
studied the problems over a three dimensional half space with various boundary conditions.
In particular, the results in \cite{AMR,GHR2,SS1} clarified that
the thickness of the boundary layer is of order of the Debye length.

It is of greater interest to study the cases of nonplanar walls, 
and know how the Bohm criterion and the thickness of the boundary layer depend on the shape of walls.
For these questions, Jung--Kwon--Suzuki in \cite{JKS1} established the existence of spherical symmetry stationary solutions of the Euler--Poisson system over an annulus, and then proposed {\it a Bohm criterion for the annulus},
which essentially differs from the original Bohm criterion \eqref{Bohm1}.
Moreover, the Bohm criterion is analyzed in a perturbed half space.
Suzuki--Takayama \cite{MM2} proved that the stationary solution exists uniquely and it is time-asymptotically stable, if the Bohm criterion and a certain necessary condition hold.
From these results, it has been seen that the Bohm criterion varies according to the shape of walls.
On the other hand, Jung--Kwon--Suzuki \cite{JKS1,JKS20,JKS21} studied the quasi-neutral limit over an annulus, 
and concluded that the thickness  is the same as in the case of planer walls, i.e., of order of the Debye length.

As mentioned above, the sheath formation for the isothermal or isentropic flow has been well studied.
Recently, Duan--Yin--Zhu \cite{DYZ1} gave the first result which investigates the sheath formation by the {\it nonisentropic} Euler--Poisson system.
They showed the unique existence and asymptotic stability of the monotone stationary solutions over a half line assuming the Bohm criterion \eqref{Bohm1}.
It would be worth extending this result to that for nonplanar walls. 
In this paper, we show the unique existence and asymptotic stability of the stationary solutions of \eqref{eq0} over a perturbed half space.

For readers' convenience,
we mention briefly researches from other perspectives on the sheath.
The interesting subject is to analyze the interface or transition between the plasma and sheath.
Riemann--Daube \cite{RD1} proposed a certain hydrodynamic model describing
the dynamics of an interface between the plasma and sheath over a half line.
Feldman--Ha--Slemrod \cite{FHS1} derived the generalized model
for the nonplanar wall cases, and studied the time-local solvability.
Liu--Slemrod in \cite{LS1} derived a certain KdV equation over a half line, which describes the transition. 

Of course, it is also intriguing to study the sheath formation by using kinetic models such as the Vlasov--Poisson system.
Badsi--Campos Pinto--Despr\'es \cite{BCD1} showed the unique existence of the stationary solution of the Vlasov--Poisson system in a bounded interval assuming {\it the kinetic Bohm criterion} proposed in Boyd--Thompson \cite{BT1}. 
Badsi \cite{Ba} also established the linear stability of the stationary solution.
Furthermore, the papers \cite{MM2,STZ1} studied the unique existence and nonlinear stability of the stationary solution in a half line under the kinetic Bohm criterion.

\medskip

\noindent
{\bf Notation.} 
The notation $\langle u,v \rangle$ means the inner product of $u,v \in \mathbb R^n$ for $n \in \mathbb N$.
We use $c$ and $C$ to denote generic positive constants.
Let us also denote a generic positive constant depending additionally
on other parameters $\alpha$, $\beta$, $\ldots$
by $C_{\alpha,\, \beta,\, \ldots}$. 
For a nonnegative integer $k$, ${\cal {B}}^k (\Sigma)$ stands for
the space of functions whose derivatives up to  
$k$-th order are continuous and bounded over $\Sigma$.
Furthermore, ${\cal {B}}^\infty (\Sigma)$ is defined by
$\cap_{k=1}^\infty {\cal {B}}^k (\Sigma)$.
For $ 1 \leq p \leq \infty$ and a nonnegative integer $k$, 
$L^p(\Omega)$ is the Lebesgue space;
$W^{k,p}(\Omega)$ is the $k$-th order Sobolev space in the $L^p$ sense;
$H^k(\Omega)$ is the $k$-th order Sobolev space in the $L^2$ sense, 
equipped with the norm $\Vert\cdot\Vert_k$.
We note $H^0 = L^2$, 
$\Vert\cdot\Vert:=\Vert\cdot\Vert_0$, and $H^\infty:=\cap_{k=1}^\infty H^k$.
We also define the exponential weighted Sobolev space $H^k_{\alpha}(\Omega)$ 
for $\alpha > 0$ by 
\begin{gather*}
H^k_{\alpha} (\Omega):=\left\{ f \in H^k(\Omega) \,\left| \, 
\|f\|_{k,\alpha}^2=\sum_{j=0}^k\int_{\Omega} e^{\alpha x_1} |\nabla^j f|^2\, dx < \infty \right\}\right..
\end{gather*}
Note that there exist $c$ and $C$ independent of $\alpha$ such that
\begin{equation}\label{eqiv0}
c\|f\|_{k,\alpha}\leq\|e^{\alpha x_1/2}f\|_k\leq C\|f\|_{k,\alpha}
\quad \text{for $f \in H^k_{\alpha} (\Omega)$ and $\alpha \in (0,1]$}.
\end{equation}
The notation $C^k([0,T];{\mathcal H})$ means 
the space of $k$-times continuously
differentiable functions on the interval $[0,T]$ with values in 
some Hilbert space ${\mathcal H}$.

\section{Main results}\label{S3/2}
Before mentioning our main results, we review a result in \cite{DYZ1}
which showed the unique existence of stationary solutions
over a half line $\mathbb R_+:=\{x_1>0\}$.
The stationary solution $(\tilde{\rho},\tilde{u},\tilde{\theta},\tilde{\phi})(x_1)$ solves the system
\begin{subequations}\label{sp0}
\begin{gather}
   (\tilde{\rho}\tilde{u})'=0,
   \label{sp1}\\
   m{\tilde{u}}\tilde{u}'+R\tilde{\theta}(\log \tilde{\rho})'+R\tilde{\theta}'={\tilde{\phi}}',
   \label{sp2}\\
    \tilde{u}\tilde{\theta}'+(\gamma-1)\tilde{\theta}\tilde{u}'=0,
   \label{sp4}\\
   \tilde{\phi}''=\tilde{\rho}-e^{-\tilde{\phi}}, \quad x_{1}>0
   \label{sp3}
   \end{gather}
with the conditions
\begin{gather}
  \inf_{x_1\in\mathbb R_{+}}\tilde{\rho}(x_1)>0, \quad  \inf_{x_1\in\mathbb R_{+}}\tilde{\theta}(x_1)>0, \\
  \lim_{x_1\rightarrow \infty}(\tilde{\rho},\tilde{u},\tilde{\theta},\tilde{\phi})(x_1)=(1,u_+,\theta_{+},0), \quad
  \tilde{\phi}(0)=\phi_b.
\end{gather}
\end{subequations}
Under the Bohm criterion \eqref{Bohm1},
the unique existence of stationary solutions $(\tilde{\rho},\tilde{u},\tilde{\theta},\tilde{\phi})$ was established
as in the next lemma.
\begin{lem}[\cite{DYZ1}]\label{1.1}
Let $u_+$ and $\theta_{+}$ satisfy \eqref{Bohm1}. There exist a constant $\delta>0$ such that
if $|\phi_b| < \delta$, then the problem \eqref{sp0}
has a unique monotone solution $(\tilde{\rho} ,\tilde{u},\tilde{\theta}, \tilde{\phi}) \in 
{\cal B}^{\infty}(\overline{\mathbb{R}_{+}})$. Moreover, it satisfies
\begin{equation}\label{ses1}
|\partial_{x_1}^j(\tilde{\rho}-1)|
+|\partial_{x_1}^j(\tilde{u}-u_{+})|
+|\partial_{x_1}^j(\tilde{\theta}-\theta_{+})|
+|\partial_{x_1}^j\tilde{\phi}|
 \leq C|\phi_b| e^{-\alpha x_{1}}
\quad \text{for} \quad 
j=0,1,2,\cdots, 
\end{equation}
where $\alpha<1$ and $C$ are positive constants independent of $\phi_b$.
\end{lem}

From now on we discuss our main results.
We first show the unique existence of stationary solutions 
$(\rho^s,\bm{u}^s,\theta^{s},\phi^s)=(\rho^s,u^s_1,u^s_2,u^s_3,\theta^{s},\phi^s)$
over the perturbed half space $\Omega$ by regarding 
$(\rho^s,u^s_1,u^s_2,u^s_3,\theta^{s},\phi^s)(x)$ 
as a perturbation of $(\tilde{\rho},\tilde{u},0,0,\tilde{\theta},\tilde{\phi})(\tilde{M}(x))$, where
\begin{equation}\label{tM1}
\tilde{M}(x):=x_1-M(x'). 
\end{equation}
The stationary solutions satisfy the equations
\begin{subequations}\label{seq0}
   \begin{gather}
     \nabla \cdot(\rho^s \bm{u}^s)=0,
   \label{seq1}\\
   m\left( \bm{u}^{s} \cdot \nabla \right) \bm{u}^{s}
    + R \theta^{s} \nabla (\log \rho^{s}) + R \nabla \theta^{s} = \nabla \phi^{s},
   \label{seq2}\\
    \bm{u}^{s} \cdot \nabla \theta^{s} 
    + (\gamma-1) \theta^{s} \nabla \cdot \bm{u}^{s} = 0,
   \label{seq4}\\
   \Delta \phi^{s}=\rho^{s}-e^{-\phi^{s}}, \quad x \in \Omega
   \label{seq3}
    \end{gather}
and the conditions
   \begin{gather}
  \inf_{x\in \Omega}\rho^s(x)>0, \quad   \inf_{x\in \Omega}\theta^s(x)>0,
  \label{spo1}
  \\
  \lim_{x_1\to\infty}(\rho^s,u^s_1,u^s_2,u^s_3,\theta,\phi)(t,x_1,x')=(1,u_+,0,0,\theta_{+},0),
  \label{sbc1}
  \\
  \phi^s(t,M(x'),x')=\phi_b \quad \text{for} \
  x'\in \mathbb R^2.
  \label{sbc2}
 \end{gather}
 \end{subequations}
The existence result is summarized in the next theorem.

\begin{thm}\label{1.2}
Let $r \geq 3$ be an integer. Suppose that $u_+$ and $\theta_{+}$ satisfy \eqref{Bohm1} and \eqref{asp1'}.
There exist positive constants $\beta \leq \alpha/2$,
where $\alpha$ is defined in Lemma \ref{1.1},
and $\delta$ such that if $|\phi_b| \leq \delta$, 
then the stationary problem \eqref{seq0} 
has a unique solution $(\rho^s,\bm{u}^s,\theta^{s},\phi^s)$ that satisfies
\begin{gather*}
(\rho^s,u^s_1,u^s_2,u^s_3,\theta^{s},\phi^s)
-(\tilde{\rho}\circ\tilde{M},\tilde{u}\circ\tilde{M},0,0,\tilde{\theta}\circ\tilde{M},\tilde{\phi}\circ\tilde{M})
\in [H^r_\beta(\Omega)]^{5}\times H^{r+1}_\beta(\Omega),
\\
\|(\rho^s-\tilde{\rho}\circ\tilde{M},u^s_1-\tilde{u}\circ\tilde{M},u^s_2,u^s_3,\theta^s-\tilde{\theta}\circ\tilde{M})\|_{r,\beta}^2
+\|\phi^s-\tilde{\phi}\circ\tilde{M}\|_{r+1,\beta}^2 \leq C |\phi_b|,
\end{gather*}
where $C$ is a positive constant independent of $\phi_b$.
\end{thm}

We also show the stability of stationary solutions
in the exponential weighted Sobolev spaces.
The paper \cite{DYZ1} pointed out that
system \eqref{eq1}--\eqref{eq3} itself does not have
the dissipative effect in the usual function space,
but there appears the effect in the weighted space.
Therefore, we employ the weighted space.

\begin{thm}\label{1.3}
Suppose that  $u_+$ and $\theta_{+}$ satisfy \eqref{Bohm1} and \eqref{asp1'}.
There exist positive constants $\beta \leq \alpha/2$,
where $\alpha$ is defined in Lemma \ref{1.1},
and $\delta$ 
such that if $\|(\rho_0-\rho^s,\bm{u}_0-\bm{u}^s,\theta_0-\theta^s)\|_{3,\beta}+|\phi_b| \leq \delta$, 
then the initial--boundary value problem \eqref{eq0} has a unique time-global solution
$(\rho,\bm{u},\theta,\phi)$ with \eqref{po1} and \eqref{super1'} in the following space:
\[
(\rho-\rho^s,\bm{u}-\bm{u}^s,\theta-\theta^s,\phi-\phi^s) \in
\left[\bigcap_{i=0}^1 C^i([0,T];H^{3-i}_\beta(\Omega))\right]^5
\times C([0,T];H^{5}_\beta(\Omega)).
\]
Moreover, there holds that
\begin{equation*}
\sup_{x \in \Omega}|(\rho-\rho^s,\bm{u}-\bm{u}^s,\theta-\theta^s,\phi-\phi^s)(t,x)| \leq Ce^{-\lambda t}
\quad \text{for $t \in [0,\infty)$},
\end{equation*}
where $C$ and $\lambda$ are positive constants independent of $\phi_b$ and $t$.
\end{thm}

\begin{rem} {\rm
What most interests us in Theorems \ref{1.2} and \ref{1.3} is that
the Bohm criterion \eqref{Bohm1} with the supersonic outflow condition \eqref{asp1'} 
still guarantees the formation of sheaths even for the perturbed half space.
The thing to note here is that \eqref{asp1'} is a necessary condition
for the well-posedness of the problem \eqref{eq0}.
Furthermore, it is worth pointing out that we do not require 
any smallness assumptions for the function $M$ 
representing the boundary of the domain $\Omega$. 
}
\end{rem}

The main difficulty of the proof is that the stationary problem
is still given by a boundary value problem to a hyperbolic--elliptic system,
although the problem over a half space or an annulus can be reduced to 
a system of ordinary differential equations.
One may first try to linearize the stationary Euler--Poisson system so that
the hyperbolic and elliptic parts are decoupled, 
and then apply an inductive scheme to solve the nonlinear problem.
However, this approach does not work well under the physically relevent situation.
(It may work if we assume that the end state $u_{+}$ is much greater than the Bohm criterion \eqref{Bohm1}.)
For the same reason, the contraction mapping principle is also not useful for our situation.

To resolve this difficulty, we borrow the approach used in \cite{MM1,Va83}.
We first show the time-global solvability of the problem \eqref{eq0}, and 
then construct stationary solutions by using the time-global solutions.
Indeed, the approach works well for our situation as follows.
In the case $\phi_b=0$, we see that the end state $(1,u_{+},0,0,,\theta_{+})$ solves the stationary problem \eqref{seq0}.
Then from the stability analysis in \cite{DYZ1} studying the nonisentropic flow, it is highly expected that the solution $(\rho,\bm{u},\theta)$ of \eqref{eq0} exists globally in time, and converges to the end state exponentially fast 
in the exponential weighted Sobolev space as $t$ tends to infinity.
On the other hand, for the case $\phi_b\neq 0$, after suitable reformulation, 
all effects coming from $\phi_b\neq 0$ are represented by inhomogeneous terms in the system.
Therefore, the dissipative structure for the end state enables us to estimate
the solutions with $\phi_b\neq 0$ by the initial data and inhomogeneous terms independent of time $t$.
For the construction of the stationary solution of \eqref{seq0}, 
we define a sequence by the time-global solution
whose time $t$ is translated to $t+kT^*$ for any $T^*>0$ and $k \in \mathbb N$, and then show that 
this sequence converges a time-periodic solution with a period $T^*$ as $k$ tends to infinity. 
Using the arbitrary of the period $T^{*}$, we conclude that
this time-periodic solution is independent of $t$, and thus obtain the stationary solution.

This paper is organized as follows.
In Section \ref{S2}, we start from rewriting the initial--boundary value problem \eqref{eq0}
by introducing a perturbation from the stationary solution
over the half space.
Section \ref{S4} is devoted to showing the time-global solvability of the rewritten problem
in the exponential weighted Sobolev space.
We construct stationary solutions in Section \ref{S5} 
by using the time-global solutions established above.
The stability of stationary solutions is also shown in the same weighted space.

\section{Reformulation and Preliminary}\label{S2}
\subsection{Reformulation}
For mathematical convenience, we begin by reformulating the initial--boundary value problem~\eqref{eq0}.
Let us introduce new functions 
\begin{align*}
v(t,x)&:=\log\rho(t,x), & \tilde{v}(x_{1})&:=\log\tilde{\rho}(x_{1}),
\\
{V}(t,x)&:={}^{t}({v},{\bm{u}}, {\theta})(t,x), & \tilde{V}(x_1)&:={}^{t}(\tilde{v},\tilde{\bm{u}}, \tilde{\theta})(x_1):={}^{t}(\log\tilde{\rho},\tilde{u},0,0, \tilde{\theta})(x_1)
\end{align*}
and perturbations
\begin{align*}
\Psi(t,x)&={}^t(\psi,\bm{\eta}, \zeta)(t,x)={}^t(\psi,\eta_1,\eta_2,\eta_3,\zeta)(t,x)
:=V(t,x)-\tilde{V}(\tilde{M}(x)),
\\
\sigma(t,x)&:=\phi(t,x)-\tilde{\phi}(\tilde{M}(x)),
\end{align*}
where $\tilde{M}(x)$ is defined in \eqref{tM1}.
Then, from \eqref{eq0} and \eqref{sp0}, we have the reformulated problem for $(\Psi,\sigma)$:
\begin{subequations}\label{re0}
\begin{gather}
A^0[{V}]\partial_{t}\Psi+
\sum_{j=1}^3A^j[{V}]\partial_{x_j}\Psi
=\begin{bmatrix}
0 \\ \nabla \sigma \\0
\end{bmatrix}
+B[V,\tilde{V}',\nabla M]\Psi
+\begin{bmatrix}
0 \\ \bm{h}[\tilde{V},\tilde{V}',\nabla M]\\0
\end{bmatrix},
\label{re1}
\\
\Delta \sigma-\sigma={\psi}+g_0[\psi,\tilde{v}]
+g_1[\sigma,\tilde{\phi}]+g_2[\tilde{\phi}',\nabla M],
\label{re2}
\\
 \lim_{|x|\to\infty}(\Psi,\sigma)(t,x)=0, 
\label{re4}
\\
 \sigma(t,M(x'),x')=0,
\label{re5}
\\
\Psi(0,x)=\Psi_0(x):=
{}^t(\log\rho_0,\bm{u}_0,\theta_{0})(x)-\tilde{V}(\tilde{M}(x)).
\label{re3}
\end{gather}
\end{subequations}
Here the $5 \times 5$ symmetric matrices $A^j$, $5 \times 5$ matrix $B$,
and $3 \times 1$ matrix $\bm{h}$ are defined as
\begin{align*}
&A^0[{V}]\!:=\!
\begin{bmatrix}
R\theta & 0 & 0 & 0 & 0
\\
0 & m & 0 & 0 & 0
\\
0 & 0 & m & 0 & 0
\\
0 & 0 & 0 & m & 0
\\
0 & 0 & 0 & 0 & \frac{R}{(\gamma-1)\theta}
\end{bmatrix}\!\!,
\quad
A^1[{V}]\!:=\!
\begin{bmatrix}
R{\theta}u_{1} & R\theta & 0 & 0 & 0
\\
R\theta & mu_1 & 0 & 0 & R
\\
0 & 0 & mu_1 & 0 & 0
\\
0 & 0 & 0 & mu_1 & 0
\\
0 & R & 0 & 0   & \frac{Ru_1}{(\gamma-1)\theta}
\end{bmatrix}\!\!,
\\
&A^2[{V}]\!:=\!
\begin{bmatrix}
R\theta u_2 & 0 & R\theta & 0 & 0
\\
0 & mu_2 & 0 & 0 & 0
\\
R\theta & 0 & mu_2 & 0 & R
\\
0 & 0 & 0  & mu_2 & 0
\\
0 & 0 & R & 0 & \frac{Ru_2}{(\gamma-1)\theta}
\end{bmatrix}\!\!,
\quad
A^3[{V}]\!:=\!
\begin{bmatrix}
R\theta u_3 & 0 & 0 & R\theta & 0
\\
0 & mu_3 & 0 & 0 & 0
\\
0 & 0  & mu_3 & 0 & 0
\\
R\theta & 0 & 0 & mu_3 & R
\\
0 & 0 & 0 & R & \frac{Ru_3}{(\gamma-1)\theta}
\end{bmatrix}\!\!,
\\
&B[{V},\tilde{V}',\nabla M]\!:=\! A^{0}[V]\tilde{B}[\tilde{V}',\nabla M],
\quad
\tilde{B}[\tilde{V}',\nabla M]\!:=\!
\begin{bmatrix}
0 & - \tilde{v}' &  \tilde{v}' M_{x_{2}} & \tilde{v}' M_{x_{3}} & 0
\\
0 & -\tilde{u}' & \tilde{u}' M_{x_{2}} & \tilde{u}' M_{x_{3}} & -\frac{R}{m}\tilde{v}'
\\[3pt]
0 & 0 & 0 & 0 & \frac{R}{m}\tilde{v}'  M_{x_{2}} 
\\[3pt]
0 & 0 & 0 & 0 & \frac{R}{m}\tilde{v}' M_{x_{3}} 
\\[3pt]
0 & -\tilde{\theta}' & \tilde{\theta}' M_{x_{2}}  & \tilde{\theta}' M_{x_{3}} & -(\gamma-1)\tilde{u}'
\end{bmatrix}\!\!,
\\
&\bm{h}[\tilde{V},\tilde{V}',\nabla M]\!:=\!
\begin{bmatrix}
0
\\
-m\tilde{u}\tilde{u}'  M_{x_{2}}
\\
-m\tilde{u}\tilde{u}' M_{x_{3}}
\end{bmatrix}\!\!.
\end{align*}
The scalar values $g_0$, $g_1$, and $g_2$ are defined as
\begin{align*}
g_0[\psi,\tilde{v}]
&:=(e^{\tilde{v}}-1){\psi}
+e^{\tilde{v}}\left(e^{\psi}-1- {\psi}\right),
\\
g_1[\sigma,\tilde{\phi}]
&:=(e^{-\tilde{\phi}}-1)\sigma-e^{-\tilde{\phi}}(e^{-\sigma}-1+\sigma),
\\
g_2[\tilde{\phi}',\nabla M]
&:=\sum_{i=2}^3 (-\tilde{\phi}''( M_{x_{i} })^2+\tilde{\phi}' M_{x_{i} x_{i} }).
\end{align*}

We also rewrite the conditions \eqref{po1} and \eqref{super1'} for this system.
The positivity of the density always holds thanks to $\rho=e^{\tilde{v}+\psi}$.
On the other hand, the positivity of the temperature is written as
\begin{gather}\label{pt1}
\inf_{x\in \Omega}(\tilde{\theta}+\zeta)(t,x)>0.
\end{gather}
Furthermore, it is straightforward to check that \eqref{super1'} is equivalent to 
\begin{gather}
\inf_{x \in \partial\Omega, \, \Phi \in \mathbb R^5, \, |\Phi|=1 } 
\left\langle \sum_{j=1}^3 n_j(x') A^j[V(t,x)]\Phi, 
\Phi \right \rangle>0.
\label{super1}
\end{gather}

It suffices to show Theorems \ref{2.1} and \ref{2.2} below 
for the completion of the proof of Theorems \ref{1.2} and \ref{1.3}, respectively.
\begin{thm}\label{2.1}
Let $r\geq 3$ be an integer. Suppose that $u_+$ and $\theta_{+}$ satisfy \eqref{Bohm1} and \eqref{asp1'}.
There exist positive constants $\beta \leq \alpha/2$,
where $\alpha$ is defined in Theorem \ref{1.1},
and $\delta$ such that if $|\phi_b| \leq \delta$, 
then the associated stationary problem of \eqref{re0} 
has a solution $(\Psi^s,\sigma^s) \in [H^r_\beta(\Omega)]^5\times H^{r+1}_\beta(\Omega)$ that satisfies \eqref{pt1}, \eqref{super1}, and
\begin{equation*}
\|\Psi^s\|_{r,\beta}^2 +\|\sigma^s\|_{r+1,\beta}^2 \leq C |\phi_b|,
\end{equation*}
where $C$ is a positive constant independent of $\phi_b$.
\end{thm}
\begin{thm}\label{2.2}
Suppose that $u_+$ and $\theta_{+}$ satisfy \eqref{Bohm1} and \eqref{asp1'}.
There exist positive constants $\beta \leq \alpha/2$,
where $\alpha$ is defined in Theorem \ref{1.1},
and $\delta$ 
such that if $\|\Psi_0\|_{3,\beta}+|\phi_b| \leq \delta$, 
then the initial--boundary value problem \eqref{re0} has a unique time-global solution
$(\Psi,\sigma) \in [\bigcap_{i=0}^1 C^i([0,T]$ $;H^{3-i}_\beta(\Omega))]^5
\times C([0,T];H^{5}_\beta(\Omega))$ 
with \eqref{pt1} and \eqref{super1}. Moreover, there holds that
\begin{equation}\label{decay1}
\sup_{x \in \Omega}|(\Psi-\Psi^s,\sigma-\sigma^s)(t,x)| \leq Ce^{-\lambda t}
\quad \text{for $t \in [0,\infty)$,}
\end{equation}
where $C$ and $\lambda$ are positive constants independent of $\phi_b$ and $t$.
\end{thm}

\subsection{Preliminary}

We use general inequalities in \cite[Lemma B.1]{MM1} in proving Theorems \ref{2.1} and \ref{2.2}.

\begin{lem}[\cite{MM1}]
Let $l=0,1,2,\cdots$ and $\beta \in [0,1]$.
Suppose that $A \in {\cal B}^\infty(B(0,r))$, 
$A(0)=0$, and $\tilde{A} \in {\cal B}^{l+1}(\overline{\Omega})$,
where $B(0,r)\subset \mathbb R^n$ denotes a ball
of center $O$ and radius $r \in (0,1]$.
If $f\in L^\infty(\Omega) \cap H^l(\Omega)$,
$g\in H^l_\beta(\Omega)$, and $e^{\beta x_1/2}g \in L^\infty(\Omega)$, there holds that
\begin{gather*}
\|fg\|_{l,\beta}
\leq C(\|f\|_{L^\infty}\|g\|_{l,\beta}
+\|f\|_{l}\|e^{\beta x_1/2}g\|_{L^\infty}),
\\
\|A(f)\|_l \leq C\|f\|_l 
\quad \text{if $\|f\|_{L^\infty} \leq r/2$}.
\end{gather*}
If $f,\nabla f\in L^\infty(\Omega) \cap H^{l}(\Omega)$,
$g\in H^l_\beta(\Omega)$, and $e^{\beta x_1/2}g \in L^\infty(\Omega)$,
the following inequalities on the commutator $[\nabla^l,\, \cdot\,]$ hold:
\begin{gather}
\|[\nabla^{l+1},f] g \|_{0,\beta} \leq 
C(\|\nabla f\|_{L^\infty}\|g\|_{l,\beta}
+\|\nabla f\|_{l}\|e^{\beta x_1/2}g\|_{L^\infty}),
\label{A2} \\
\|[\nabla^{l+1},\tilde{A}] g \|_{0,\beta} 
\leq C \left(\sum_{i=1}^{l+1}\|\nabla^{i} \tilde{A}\|_{L^\infty} \right)
\|g\|_{l,\beta}.
\label{A4} 
\end{gather}
Here $C$ is a positive constant independent of $f$, $g$, and $\beta$. 
\end{lem}

\section{Time-global solvability}\label{S4}
This section deals with the time-global solvability of the initial--boundary value problem \eqref{re0} 
for small initial data $\Psi_0$ and boundary data $\phi_b$.

\begin{thm}\label{4.1}
Let $r \geq 3$ be an integer. 
Suppose that  $u_+$ and $\theta_+$ satisfy \eqref{Bohm1} and \eqref{asp1'}.
There exist positive constants $\beta \leq \alpha/2$,
where $\alpha$ is defined in Theorem \ref{1.1},
and $\delta$ depending on $\beta$ such that if 
$\|\Psi_0\|_{r,\beta}+|\phi_b| \leq \delta$, 
then the initial--boundary value problem \eqref{re0} has a unique time-global solution
$(\Psi,\sigma) \in [\bigcap_{i=0}^1 C^i([0,\infty);H^{r-i}_\beta(\Omega))]^5
\times C([0,\infty);H^{r+2}_\beta(\Omega))$
with \eqref{pt1} and \eqref{super1}. Moreover, there holds that
\begin{equation}\label{apes1}
\sup_{t\in[0,\infty)}\left(
\|\Psi(t)\|_{r,\beta}^2+\|\Psi_{t}(t)\|_{r-1,\beta}^2
+\|\sigma(t)\|_{r+2,\beta}^2\right)
\leq C(\|\Psi_0\|_{r,\beta}^2+|\phi_b|),
\end{equation}
where $C$ is a positive constant depending on $\beta$ but
independent of $\Phi_{0}$ and $\phi_b$.
\end{thm}

The time-global solution $(\Phi,\sigma)$ with \eqref{apes1} 
can be constructed by a standard continuation argument 
using the time-local solvability in Lemma \ref{3.1} and 
the a priori estimate in Proposition \ref{4.2} below.
For notational convenience, we use
\[
 N_{r,\alpha}(T):=\sup_{t\in [0,T]}\|\Psi(t)\|_{r,\alpha}.
\]

\begin{lem}\label{3.1}
Suppose that $\Psi_0$ satisfies \eqref{pt1} and \eqref{super1} as well as 
belongs $H^{r}_{\alpha/2}(\Omega)$ 
for $r \geq 3$ and $\alpha>0$ being in Theorem \ref{1.1}.
Let $\beta$ be a positive constant less than
$\alpha/2$ and $2 e^{m_*/2}$, where
\begin{gather*}
\label{beta1}
\quad
m_*:=\min \left\{
\inf_{x \in \mathbb R_+} \left(-\tilde{\phi}(x)\right), \
\inf_{x \in \Omega}(v_0+\tilde{v})(x)-1
\right\}.
\end{gather*}
Then there exist positive constants $\delta$ and $T$ such that 
if $|\phi_b| < \delta$, the problem \eqref{re0} has a unique solution 
$(\Psi,\sigma) \in [\bigcap_{i=0}^1 C^i([0,T];H^{r-i}_\beta(\Omega))]^5
\times C([0,T];H^{r+2}_\beta(\Omega))$ with \eqref{pt1} and \eqref{super1}.
\end{lem}

\begin{pro}\label{4.2}
Let $r \geq 3$ be an integer.  Suppose that  $u_+$ and $\theta_+$ satisfy \eqref{Bohm1} and \eqref{asp1'}.
Assume that $(\Psi,\sigma) \in [\bigcap_{i=0}^1 C^i([0,T]$ $;H^{r-i}_\beta(\Omega))]^5
\times C([0,T];H^{r+2}_\beta(\Omega))$
be a solution to the problem \eqref{re0} with \eqref{pt1} and \eqref{super1}.
There exist positive constants $\beta \leq \alpha/2$,
where $\alpha$ is defined in Theorem \ref{1.1},
and $\delta$ depending on $\beta$
such that if $N_{r,\beta}(T)+|\phi_b| < \delta$,
the following estimate holds:
\begin{equation}\label{apes2}
\sup_{t\in[0,T]}\left(
\|\Psi(t)\|_{r,\beta}^2+\|\Psi_{t}(t)\|_{r-1,\beta}^2
+\|\sigma(t)\|_{r+2,\beta}^2\right)
\leq C(\|\Psi_0\|_{r,\beta}^2+|\phi_b|),
\end{equation}
where $C$ is a positive constant depending on $\beta$
but independent of $\Phi_{0}$ and $\phi_b$.
\end{pro}

We omit the proof of Lemma \ref{3.1}, since it can be proved in much the same way as Lemma 3.1 in \cite{M.S.1}.
In the remainder of this section, we prove only Proposition \ref{4.2}.
To this end, we follow the approach used in \cite{MM1} which studied the isothermal flow over the perturbed half space $\Omega$.
More precisely, in the case $\phi_b=0$, i.e., the inhomogeneous terms $\bm{h}$ in \eqref{re1} and $g_2$ in \eqref{re2} vanish,
it is highly expected from the study \cite{DYZ1} over a half line that $(\Psi,\sigma)$ exists globally in time and decays exponentially fast 
in the exponential weighted Sobolev space as $t$ tends to infinity.
This dissipative structure also enables us to prove that
the $H^{r}_\beta$-norm of solutions with $\phi_b\neq 0$ is bounded by those of $\Psi_{0}$, $\bm{h}$, and $g_2$.
The estimate is exactly \eqref{apes2}.

\subsection{Elliptic estimates}\label{S4.1}
This subsection provides
\footnote{We remark that all constants $C$ in subsection \ref{S4.1} are independent of $\beta$.}{estimates}
of $\sigma$ solving the elliptic equation \eqref{re2}. 
First we derive the lower and upper bounds of $\sigma$.
\begin{lem}\label{ell0}
Under the same assumption as in Proposition \ref{4.2}, there hold that
\begin{gather*}
\sup_{x\in\Omega}(\sigma+\tilde{\phi})(t,x) \leq M_1,
\quad M_1:=\max\left\{\sup_{x\in\Omega}|\tilde{\phi}(\tilde{M}(x))|, \
-\inf_{x\in\Omega}{\tilde{v}(\tilde{M}(x))}+1\right\},
\\
\inf_{x\in\Omega}(\sigma+\tilde{\phi})(t,x) \geq -M_2,
\quad M_2:=\max\left\{\sup_{x\in\Omega}|\tilde{\phi}(\tilde{M}(x))|, \
\sup_{x\in\Omega}\tilde{v}(\tilde{M}(x))+1\right\},
\\
\sup_{x\in\Omega}|\sigma(t,x)| \leq C (N_{r,\beta}(T)+|\phi_b|),
\end{gather*}
where $C$ is a positive constant independent of $\beta$, $\phi_b$, and $t$.
\end{lem}
\begin{proof}
We can show this lemma in much the same way as the proof of \cite[Lemma 4.4]{MM1}.
Indeed, it follows from taking $K=1$ in the proof of \cite[Lemma 4.4]{MM1}.
\end{proof}

Next we also obtain the estimate of the $H^k_\beta$-norm of $\sigma$.

\begin{lem}\label{ell1}
Under the same assumption as in Proposition \ref{4.2}, there hold that
\begin{gather}
\|\sigma(t)\|_{1,\beta}^2 
\leq \{1+D\beta^2+C(N_{r,\beta}(T)+|\phi_b|)\} \|\psi(t)\|_{0,\beta}^2+C|\phi_b|,
\label{ellineq1}\\
\|\sigma(t)\|_{l+2,\beta}^2 
\leq C(\|\psi(t)\|_{l,\beta}^2+|\phi_b|^2) \quad \text{for $l=0,1,\ldots,r$,}
\label{ellineq4}
\end{gather}
where $C$ and $D$ are positive constants independent of $\beta$, $\phi_b$, and $t$.
\end{lem}
\begin{proof}
This lemma can be shown in much the same way as the proof of \cite[Lemma 4.5]{MM1}.
\end{proof}

\subsection{Basic estimate}\label{S4.2}
This subsection is devoted to deriving an estimate of $L^2$-norm
of $\Psi$ solving the hyperbolic system \eqref{re1}.
Only in this subsection, we must be careful to investigate the dependence of $\beta$
in order to take it suitably small.

We begin by deriving several equalities.
Taking the inner product of \eqref{re1} with the vector $2\Psi$ gives
\begin{gather}
(\langle A^0[{V}]\Psi,\Psi \rangle)_{t}
+\sum_{j=1}^3\left(\langle A^j[{V}]\Psi,\Psi \rangle
-2\sigma\eta_{j}\right)_{x_{j}}
=-2\sigma(\nabla\cdot\bm{\eta})+\mathcal{R}_1,
\label{bes11}
\end{gather}
where
\begin{gather}
\mathcal{R}_1:=\langle \{ ( A^0[V] )_{t} \}\Psi,\Psi \rangle+
\sum_{j=1}^3 \langle \{ (A^j[V])_{x_{j}} \}\Psi,\Psi \rangle
+2\langle B\Psi,\Psi\rangle+2\bm{h}\cdot\bm{\eta}.
\notag
\end{gather}

Applying $\partial_{x_i}$ to the $(1+i)$th component of \eqref{re1} 
and sum up them, 
we have the following equation for $\nabla \cdot \bm{\eta}$:
\begin{gather*}
\begin{aligned}
m(\nabla \cdot  \bm{\eta})_t+R \theta \Delta \psi
+m\sum_{j=1}^3 u_j (\nabla \cdot  \bm{\eta})_{x_j}+R\Delta\zeta
=\Delta \sigma+\nabla \cdot \bm{h}
+I_1,
\end{aligned}
\\
\begin{aligned}
I_1:=-R\sum_{i=1}^3 \theta_{x_i} \psi_{x_i}-m\sum_{i,j=1}^3 u_{jx_i} \eta_{i x_j}
+\sum_{i=1}^3\{m(\tilde{B}\Psi)_{1+i}\}_{x_i}, 
\end{aligned}
\notag
\end{gather*}
where $(\tilde{B}\Psi)_{1+j}$ means the $(1+j)$th component of the vector $\tilde{B}\Psi$.
Multiplying this by $2\nabla \cdot  \bm{\eta}$ leads to 
\begin{align}\label{bes33}
&m\{(\nabla \cdot  \bm{\eta})^2 \}_t
+2R\theta(\Delta \psi) (\nabla \cdot  \bm{\eta})
+m\sum_{j=1}^3 \{ u_j (\nabla \cdot  \bm{\eta})^2\}_{x_j}+2R\Delta\zeta(\nabla \cdot  \bm{\eta})
\notag\\
&=2\Delta \sigma(\nabla \cdot  \bm{\eta})+{\mathcal R}_2,
\end{align}
where
\begin{equation*}
{\mathcal R}_2:=m\sum_{j=1}^3 (u_{j})_{x_j}(\nabla \cdot  \bm{\eta})^2
+2(\nabla \cdot \bm{h})(\nabla \cdot  \bm{\eta})
+2I_1(\nabla \cdot  \bm{\eta}).
\end{equation*}

We multiply the first component of \eqref{re1} by $(R\theta)^{-1}$, 
apply the operator $\nabla$, and multiply the resultant by $R\theta$ to obtain the following equations for $\nabla \psi$:
\begin{gather*}
R \theta (\nabla \psi)_t+R \theta \sum_{j=1}^3 u_{j}(\nabla \psi)_{x_j}
+R \theta \nabla(\nabla \cdot  \bm{\eta})
=I_2
:= -R \theta \nabla \bm{u} \nabla \psi
+R\theta \nabla(\tilde{B}\Psi)_1, 
\end{gather*}
where $(\tilde{B}\Psi)_1$ means the 1st component of the vector $\tilde{B}\Psi$.
Taking inner product of this with $2\nabla \psi$ gives
\begin{equation}\label{bes55}
\{R \theta |\nabla \psi|^2 \}_t
+\sum_{j=1}^3\left\{R \theta u_{j} |\nabla \psi|^2 \right\}_{x_j}
+2\nabla\cdot\{R\theta(\nabla \cdot  \bm{\eta}) \nabla \psi\}
-2R\theta(\Delta \psi)(\nabla \cdot  \bm{\eta})
={\mathcal R}_3,
\end{equation}
where
\begin{equation*}
{\mathcal R}_3:=R\theta_t|\nabla \psi|^2+
\sum_{j=1}^3 \left\{R \theta u_{j} \right\}_{x_j}|\nabla \psi|^2
+2R(\nabla \cdot  \bm{\eta}) (\nabla \theta \cdot \nabla \psi)
+2I_2\cdot \nabla \psi.
\end{equation*}


We multiply the fifth component of \eqref{re1} by $(\gamma-1)R^{-1}\theta$, 
apply the operator $\nabla$, and multiply the resultant by $(\gamma-1)^{-1}R\theta^{-1}$ to obtain the following equations for $\nabla \zeta$:
\begin{gather*}
\frac{R}{(\gamma-1)\theta}(\nabla\zeta)_t+\frac{R}{(\gamma-1)\theta}\sum_{j=1}^3 u_{j}(\nabla\zeta)_{x_j}
+R\nabla(\nabla \cdot  \bm{\eta})
=I_{3},
\\
\begin{aligned}
I_3
&:=-\frac{R\nabla \bm{u} \nabla \zeta}{(\gamma-1)\theta }
-\frac{R(\nabla\cdot\bm{\eta})\nabla\theta}{\theta}
+\frac{R\nabla(\tilde{B}\Psi)_5}{(\gamma-1)\theta}.
\end{aligned}
\notag
\end{gather*}
Taking inner product of this with $2\nabla \zeta$ leads to 
\begin{gather}
\left\{\frac{R}{(\gamma-1)\theta }|\nabla\zeta|^2\right\}_t+\sum_{j=1}^3\left\{\frac{Ru_j}{(\gamma-1)\theta} |\nabla\zeta|^2\right\}_{x_j}
+2\nabla\cdot\{R (\nabla \cdot  \bm{\eta}) \nabla\zeta\}-2R(\nabla \cdot  \bm{\eta})\Delta\zeta
={\mathcal R}_4,
\label{bes66}
\end{gather}
where
\begin{gather*}
{\mathcal R}_4:=\left\{\frac{R}{(\gamma-1)\theta }\right\}_t|\nabla \zeta|^2+
\sum_{j=1}^3 \left\{\frac{Ru_{j}}{(\gamma-1)\theta}\right\}_{x_j}|\nabla \zeta|^2
+2I_3\cdot \nabla \zeta.
\end{gather*}


To handle the terms having $\sigma$ 
on the right hand sides of \eqref{bes11} and \eqref{bes33},
we multiply \eqref{re2} by $2\nabla\cdot\bm{\eta}$ and rewrite the result as
\begin{align}
2(\Delta \sigma-\sigma)(\nabla\cdot\bm{\eta})
&=2\psi(\nabla\cdot\bm{\eta})
+2(g_0+g_1+g_2)(\nabla\cdot\bm{\eta})
\notag\\
&=2\psi\left(-\psi_{t}
-\sum_{j=1}^3 u_{j} \psi_{x_{j}}
+(\tilde{B}\Psi)_1\right)
+2(g_0+g_1+g_2)(\nabla\cdot\bm{\eta})
\notag\\
&=-(\psi^2)_{t}
-\sum_{j=1}^3\left\{u_{j}\psi^2 \right\}_{x_{j}}
+{\mathcal R}_5,
\label{bes4}
\end{align}

\vspace{-4mm}

\[
{\mathcal R}_5:= (\nabla\cdot\bm{u}) \psi^2
+2\psi(\tilde{B}\Psi)_1+2(g_0+g_1+g_2)(\nabla\cdot\bm{\eta}),
\]
where we have also used the first component of \eqref{re1}
in deriving the second equality.

Summing up the equalities \eqref{bes11}--\eqref{bes4}, we arrive at
\begin{align}
{}&
\left(\langle  A^0[V]\Psi,\Psi\rangle+\psi^2+m(\nabla \cdot  \bm{\eta})^2+R \theta |\nabla \psi|^2+\frac{R}{(\gamma-1) \theta }|\nabla\zeta|^2\right)_t
\notag\\
&\quad +2\nabla\cdot\{R \theta (\nabla \cdot  \bm{\eta}) \nabla \psi+R(\nabla \cdot  \bm{\eta}) \nabla \zeta \}
\notag\\
&\quad +\sum_{j=1}^3\left( \langle A^j[V]\Psi,\Psi \rangle  + u_{j}\psi^2 - 2\sigma\eta_{j}
+m u_j (\nabla \cdot  \bm{\eta})^2 
+R \theta u_{j} |\nabla \psi|^2+
\frac{R}{(\gamma-1)\theta} u_{j}|\nabla\zeta|^2\right)_{x_j}
\notag\\
&
=\sum_{i=1}^5 {\mathcal R}_i.
\label{bes8}
\end{align}
Then we multiply \eqref{bes8} by $e^{\beta x_1}$, integrate it over $\Omega$, 
and use Gauss's divergence theorem with the boundary condition \eqref{re5} 
to obtain
\begin{align}
{}&
\frac{d}{dt} \int_\Omega e^{\beta x_1}
\left(\psi^2+
\langle  A^0[{V}]\Psi,\Psi\rangle+m(\nabla \cdot  \bm{\eta})^2+R\theta |\nabla \psi|^2+\frac{R}{(\gamma-1)\theta}|\nabla\zeta|^2\right)\,dx
\notag\\
&
\quad+\sum_{j=1}^3\int_{\partial\Omega} e^{\beta M(x')}
\langle n_j A^j[V]\Psi,\Psi \rangle \,ds
+\sum_{j=1}^3\int_{\partial\Omega} e^{\beta M(x')}n_ju_{j}\psi^2\,ds
\notag\\
&
\quad+\int_{\partial\Omega} e^{\beta M(x')}\left\langle\!
F[V,\bm{n}]
\begin{bmatrix}
\nabla \cdot  \bm{\eta}
\\
\nabla \psi
\\
\nabla \zeta
\end{bmatrix},
\begin{bmatrix}
\nabla \cdot  \bm{\eta}
\\
\nabla \psi
\\
\nabla \zeta
\end{bmatrix}
\!\right\rangle\,ds
\notag\\
&
\quad-\beta\int_{\Omega} e^{\beta x_1}
(\langle A^1[V]\Psi,\Psi \rangle
+u_{1}\psi^2
-2\sigma\eta_1)\,dx
-\beta\int_{\Omega} e^{\beta x_1}
\left\langle\!
F^1[{V}]
\begin{bmatrix}
\nabla \cdot  \bm{\eta}
\\
\nabla \psi
\\
\nabla \zeta
\end{bmatrix},
\begin{bmatrix}
\nabla \cdot  \bm{\eta}
\\
\nabla \psi
\\
\nabla \zeta
\end{bmatrix}
\!\right\rangle\,dx
\notag\\
&=\int_{\Omega} e^{\beta x_1} \sum_{i=1}^5 {\mathcal R}_i \,dx.
\label{bes7}
\end{align}
Here the $7 \times 7$ symmetric matrices $F$ and $F^{1}$ are defined by
\begin{align}
F[V,\bm{n}]&:=
\begin{bmatrix}
m \bm{n}\cdot\bm{u} & n_1R\theta &  n_2R\theta & n_3R\theta & n_1R &n_2R & n_3R
\\
 n_1R\theta & R\theta \bm{n}\cdot\bm{u} & 0 & 0 & 0 & 0 & 0
\\
 n_2R\theta & 0 & R\theta \bm{n}\cdot\bm{u} & 0 & 0 & 0 & 0
\\
 n_3R\theta & 0 & 0 & R\theta \bm{n}\cdot\bm{u} & 0 & 0 & 0
 \\
 n_1R & 0 & 0 & 0 & \frac{R\bm{n}\cdot\bm{u}}{(\gamma-1)\theta} & 0 & 0
\\
 n_2R & 0 & 0 & 0 & 0 & \frac{R\bm{n}\cdot\bm{u}}{(\gamma-1)\theta} & 0 
 \\
 n_3R & 0 & 0 & 0 & 0 & 0 & \frac{R\bm{n}\cdot\bm{u}}{(\gamma-1)\theta} 
\end{bmatrix},
\label{DefF}\\
F^1[V]&:=
\begin{bmatrix}
m u_1 & R\theta &  0 & 0 & R & 0 & 0
\\
 R\theta & R\theta u_1 & 0 & 0 & 0 & 0 & 0
\\
 0 & 0 & R\theta u_1 & 0 & 0 & 0 & 0
\\
 0 & 0 & 0 & R\theta u_1 & 0 & 0 & 0
 \\
R & 0 & 0 & 0 & \frac{Ru_1}{(\gamma-1)\theta} & 0 & 0
\\
 0 & 0 & 0 & 0 & 0 & \frac{Ru_1}{(\gamma-1)\theta} & 0 
 \\
 0 & 0 & 0 & 0 & 0 & 0 & \frac{Ru_1}{(\gamma-1)\theta} 
\end{bmatrix}.
\label{DefF1}
\end{align}
It follows from \eqref{super1'} and \eqref{Bohm1} that
\begin{gather}
\inf_{x \in \partial\Omega, \, \Phi \in \mathbb R^7, \, |\Phi|=1 } 
\left\langle F[V(t,x),\bm{n}(x)]\Phi, 
\Phi \right \rangle>0,
\label{F0}\\
\inf_{\Phi \in \mathbb R^7, \, |\Phi|=1 } 
\left\langle -F^{1}[V_{+}]\Phi, 
\Phi \right \rangle>0,
\label{F1}
\end{gather}
where $V_{+}:={}^{t}(1,u_{+},0,0,\theta_{+})$.

From now on we estimate the $L^2$-norm of $\Psi$.

\begin{lem}\label{4.3}
Under the same assumption as in Proposition \ref{4.2},
there holds that 
\begin{equation}\label{basic1}
\sup_{t\in[0,T]}\|\Psi(t)\|_{0,\beta}^2
\leq C\|\Psi_0\|_{1,\beta}^2
+\frac{C}{\beta}(N_{r,\beta}(T)+|\phi_b|)\sup_{t\in[0,T]}\|\nabla \Psi(t)\|_{0,\beta}^2
+\frac{C}{\beta}|\phi_b|,
\end{equation}
where $C$ is a positive constant independent of $\beta$, $\phi_b$, and $t$.
\end{lem}
\begin{proof}
To obtain \eqref{basic1}, we estimate each terms on the left hand side of \eqref{bes7} from below separately.
The second and fourth terms are nonnegative owing to \eqref{super1} and \eqref{F0}, respectively.
We also see from \eqref{ses1} and $n_1u_+>0$ that the third term is also nonnegative as follows:
\begin{align*}
&\sum_{j=1}^3\int_{\partial\Omega} e^{\beta M(x')}n_j u_{j} \psi^2\,ds
\\
&=\int_{\partial\Omega} e^{\beta M(x')}n_1 \{ u_{+} +(\tilde{u}-u_{+}) \} \psi^2\,ds
+\sum_{j=2}^3\int_{\partial\Omega} e^{\beta M(x')}n_j \eta_{j} \psi^2\,ds
\\
&\geq \int_{\partial\Omega} e^{\beta M(x')} n_1  u_+ \psi^2\,ds
-C(N_{r,\beta}(T)+|\phi_b|)\|\psi\|_{L^2(\partial \Omega)}^2 \geq 0,
\end{align*}
where we also have taking $N_{r,\beta}(T)$ and $|\phi_b|$ small enough in deriving the last inequality.
We estimate the fifth term by using \eqref{ses1} and \eqref{ellineq1} as
\begin{align*}
{}&
-\beta\int_{\Omega} e^{\beta x_1}
(\langle A^1[V]\Psi,\Psi \rangle+u_{1}\psi^2
-2\sigma\eta_1)\,dx
\\
&\geq -\beta\int_{\Omega} e^{\beta x_1}\left\{
(R\theta_++1)u_+\psi^2+2R\theta_+\psi\eta_1+2R\zeta\eta_1+m u_+|\bm{\eta}|^2+\frac{R u_+}{(\gamma-1)\theta_+}\zeta^2 -2\sigma\eta_1\right\}\,dx
\\
& \qquad
-C(N_{r,\beta}(T)+|\phi_b|)\|\Psi\|_{0,\beta}^2
\\
& \geq \beta{\mathcal D}
-\{2\sqrt{D}\beta^2+\mu+C_{\mu}(N_{r,\beta}(T)+|\phi_b|)\}\|\Psi\|_{0,\beta}^2
-C_{\mu}|\phi_b|,
\end{align*}
where $\mu$ is a positive constant to be determined later, and ${\mathcal D}$ is defined by
\begin{align*}
{\mathcal D}&:= -\int_{\Omega} e^{\beta x_1}\left\{
(R\theta_++1)u_+\psi^2+2R\theta_+\psi\eta_1+2R\zeta\eta_1+m u_+|\bm{\eta}|^2+\frac{R u_+}{(\gamma-1)\theta_+}\zeta^2\right\}\,dx
\\
&\quad -2 \|\psi\|_{0,\beta}\|\eta_1\|_{0,\beta}.
\end{align*}
By Schwarz's inequality and \eqref{Bohm1}, we see that
the term ${\mathcal D}$ is bounded from below as
\begin{align*}
{\mathcal D} 
&\geq -(R\theta_++1)u_+\|\psi\|_{0,\beta}^2
-2(R\theta_++1) \|\psi\|_{0,\beta}\|\eta_1\|_{0,\beta}
-2R \|\zeta\|_{0,\beta}\|\eta_1\|_{0,\beta}
\\
&\quad -mu_+\|\bm{\eta}\|_{0,\beta}^2 -\frac{R u_+}{(\gamma-1)\theta_+} \|\zeta\|_{0,\beta}^2
\\
&\geq d\|\Psi\|_{0,\beta}^2,
\end{align*}
where $d$ is a positive constant independent of $\beta$, $\phi_b$, and $t$.
Furthermore, one can estimate the sixth term using \eqref{F1} as follows:
\begin{align*}
{}&
-\beta\int_{\Omega} e^{\beta x_1}
\left\langle\!
F^1[V]
\begin{bmatrix}
\nabla \cdot  \bm{\eta}
\\
\nabla \psi
\\
\nabla \zeta
\end{bmatrix},
\begin{bmatrix}
\nabla \cdot  \bm{\eta}
\\
\nabla \psi
\\
\nabla \zeta
\end{bmatrix}
\!\right\rangle\,dx
\\
&\geq -\beta\int_{\Omega} e^{\beta x_1}
\left\langle\!
F^1[{V}_+]
\begin{bmatrix}
\nabla \cdot  \bm{\eta}
\\
\nabla \psi
\\
\nabla \zeta
\end{bmatrix},
\begin{bmatrix}
\nabla \cdot  \bm{\eta}
\\
\nabla \psi
\\
\nabla \zeta
\end{bmatrix}
\!\right\rangle\,dx
-C(N_{r,\beta}(T)+|\phi_b|)\|\Psi\|_{1,\beta}^2
\\
& \geq \beta d\|(\nabla \cdot \bm{\eta},\nabla \psi,\nabla\zeta)\|_{0,\beta}^2
-C(N_{r,\beta}(T)+|\phi_b|)\|\Psi\|_{1,\beta}^2.
\end{align*}
Now all terms on the left hand side except the first term 
has been estimated from below.

Next we deal with the remainder terms $\mathcal{R}_i$ on the right hand side of \eqref{bes7}.
First it is seen from \eqref{ses1} and $M \in H^\infty(\Omega)$ that
\begin{gather}
|(\bm{h},g_2,\nabla\bm{h})| \leq C|\phi_b| e^{-\alpha( x_1-M(x'))} |(\nabla M,\nabla^{2} M)|.
\label{h1}
\end{gather}
It also follows from \eqref{re1}, \eqref{ellineq4}, \eqref{h1}, and Sobolev's inequality that
\begin{align}\label{DtPsi1}
|\Psi_{t}|
= 
\left|-\sum_{j=1}^3 (A^0 )^{-1} A^j \Psi_{x_{j}}
+\begin{bmatrix}
0 \\ m^{-1} \nabla \sigma \\0
\end{bmatrix}
+ \tilde{B}\Psi
+\begin{bmatrix}
0 \\ m^{-1}\bm{h}\\0
\end{bmatrix}
\right|
 \leq C(N_{r,\beta}(T)+|\phi_b|).
\end{align}
Using \eqref{ses1}, \eqref{ellineq4}, Sobolev's and Schwarz's inequalities, and $M \in H^\infty(\Omega)$, we see that
\begin{gather*}
|\mathcal{R}_i| 
\leq C(N_{r,\beta}(T)+|\phi_b|)|(\Psi,\nabla\Psi,\sigma)|^2
+C|(\bm{h},g_2,\nabla\bm{h})||(\Psi,\nabla\Psi)|, \quad i=1,2,\cdots,5,
\end{gather*}
where $C$ is a positive constant independent of $\beta$, $\phi_b$, and $t$.
Now the right hand side of \eqref{bes7} can be estimated as follows:
\begin{gather*}
\int_{\Omega} e^{\beta x_1} \sum_{i=1}^5 {\mathcal R}_i \,dx
\leq C(N_{r,\beta}(T)+|\phi_b|)\|\Psi\|_{1,\beta}^2+C|\phi_b|,
\end{gather*}
where we have also used \eqref{ellineq4}, $\beta\leq\alpha/2$, and Schwarz's inequality.

Substituting all the above estimates into \eqref{bes7} leads to
\begin{align*}
{}&
\frac{d}{dt} \int_\Omega e^{\beta x_1}
\left(\psi^2+
\langle  A^0[{V}]\Psi,\Psi\rangle+m(\nabla \cdot  \bm{\eta})^2+R\theta |\nabla \psi|^2+\frac{R}{(\gamma-1)\theta}|\nabla\zeta|^2\right)\,dx
\\
&\quad +d\beta\|(\Psi,\nabla \psi,\nabla \cdot \bm{\eta},\nabla \zeta)\|_{0,\beta}^2
\notag\\
&\leq  2\sqrt{D}\beta^2\|\Psi\|_{0,\beta}^2
+\mu\|\Psi\|_{0,\beta}^2
+C_{\mu}(N_{r,\beta}(T)+|\phi_b|)\|\Psi\|_{1,\beta}^2+C_{\mu}|\phi_b|.
\end{align*}
To absorb the first term on the right hand side into the second term on the left hand side,
\footnote{We remark that here is only one place to choose $\beta$ suitably small and
hereafter we never change $\beta$ in the proofs of Theorems \ref{2.1} and \ref{2.2}.}{we} fix $\beta>0$ so small that
\begin{equation}\label{defbeta1}
\beta \leq \min\{\alpha/2,d(4\sqrt{D})^{-1}\}.
\end{equation}
Then taking $\mu$, $N_{r,\beta}(T)$, and $|\phi_b|$ suitably small yields
\begin{align*}
{}&
\frac{d}{dt} \int_\Omega e^{\beta x_1}
\left(\psi^2+
\langle  A^0[{V}]\Psi,\Psi\rangle+m(\nabla \cdot  \bm{\eta})^2+R\theta |\nabla \psi|^2+\frac{R}{(\gamma-1)\theta}|\nabla\zeta|^2\right)\,dx
\\
&\quad +d\beta\|(\Psi,\nabla \psi,\nabla \cdot \bm{\eta},\nabla \zeta)\|_{0,\beta}^2
\notag\\
&\leq C(N_{r,\beta}(T)+|\phi_b|)\|\nabla\Psi\|_{0,\beta}^2
+C|\phi_b|.
\end{align*}
Furthermore, multiplying this by $e^{\tilde{c} \beta t}$,
 integrating over $[0,t]$,
and taking $\tilde{c}>0$ small, we have
\begin{align*}
{}&
e^{\tilde{c} \beta t}\|(\Psi,\nabla \psi,
\nabla \cdot \bm{\eta},\nabla \zeta)(t)\|_{0,\beta}^2
+c\beta\int_0^t e^{\tilde{c} \beta \tau}
\|(\Psi,\nabla \psi,\nabla \cdot \bm{\eta},\nabla \zeta)(\tau)\|_{0,\beta}^2\,d\tau
\notag \\
&\leq C\|\Psi_0\|_{1,\beta}^2
+\int_0^t e^{\tilde{c}\beta \tau}
\left(C(N_{r,\beta}(T)+|\phi_b|)\|\nabla\Psi(\tau)\|_{0,\beta}^2
+C|\phi_b|\right)\,d\tau
\notag \\
&\leq C\|\Psi_0\|_{1,\beta}^2
+\left(C(N_{r,\beta}(T)+|\phi_b|)
\sup_{t\in[0,T]}\|\nabla\Psi(t)\|_{0,\beta}^2
+C|\phi_b|\right)
\frac{1}{\tilde{c} \beta}(e^{\tilde{c} \beta t}-1),
\end{align*}
which immediately gives \eqref{basic1}.
\end{proof}

\subsection{Higher order estimate}\label{S4.3}
In this subsection, we estimate the higher order derivatives of $\Psi$.
We multiply \eqref{re0} by $(A^{0})^{-1}$, 
apply the operator $\partial_x^{\bm{a}}$ with $|\bm{a}|=k$ for $k=1,2,\ldots,r$,
and multiply the result by $A^{0}$ to obtain
\begin{gather}
A^0[V] \partial_t(\partial_x^{\bm{a}}\Psi)+
\sum_{j=1}^3A^j[V]\partial_{x_j}(\partial_x^{\bm{a}}\Psi)=
\begin{bmatrix}
0 \\ \partial_x^{\bm{a}}\bm{h} \\ 0
\end{bmatrix}
+A^0[V]  I^{\bm{a}},
\label{hes1}\\
I^{\bm{a}}:=\sum_{j=1}^3[\partial_x^{\bm{a}},(A^0)^{-1} A^j ]\partial_{x_j}\Psi
+ \begin{bmatrix}
0 \\  m^{-1}  \nabla \partial_x^{\bm{a}} \sigma \\ 0
\end{bmatrix}
+\partial_x^{\bm{a}}(\tilde{B}\Psi),
\notag
\end{gather}
where $[\partial_x^{\bm{a}},\cdot]$ denotes a commutator.
Then we take an inner product of \eqref{hes1} with $2e^{\beta x_1}\partial_x^{\bm{a}}\Psi$,
sum up the results for $\bm{a}$ with $|\bm{a}|=k$,
integrate the resultant equality by parts over $\Omega$,
and apply Gauss's divergence theorem to obtain
\begin{align}
{}&
\frac{d}{dt} \sum_{|\bm{a}|=k}\int_\Omega e^{\beta x_1} \langle A^0[V]\partial_x^{\bm{a}}\Psi,\partial_x^{\bm{a}}\Psi \rangle \,dx
+\sum_{|\bm{a}|=k}\sum_{j=1}^3\int_{\partial\Omega} e^{\beta M(x')}
\langle n_j A^j[{V}]\partial_x^{\bm{a}}\Psi,\partial_x^{\bm{a}}\Psi \rangle\,ds
\notag\\
&\quad
-\beta\sum_{|\bm{a}|=k}\int_{\Omega} e^{\beta x_1}
\langle A^1[V]\partial_x^{\bm{a}}\Psi,\partial_x^{\bm{a}}\Psi \rangle\,dx
\notag\\
&=\sum_{|\bm{a}|=k}\int_{\Omega} e^{\beta x_1} \left(
\langle \{\partial_{t}(A^0[{V}])\}\partial_x^{\bm{a}}\Psi,\partial_x^{\bm{a}}\Psi \rangle
+\sum_{j=1}^3 \langle \{\partial_{x_j}(A^j[{V}])\}\partial_x^{\bm{a}}\Psi,\partial_x^{\bm{a}}\Psi \rangle
\right) \,dx
\notag\\
&\quad +\sum_{|\bm{a}|=k}\int_{\Omega} e^{\beta x_1} 
\left(\partial_x^{\bm{a}}\bm{h}\cdot\partial_x^{\bm{a}} \bm{\eta}
+\langle A^{0}[V]I^{\bm{a}}, \partial_x^{\bm{a}}\Psi \rangle \right) \,dx.
\label{hes4}
\end{align}

Let us estimate the higher order derivatives of $\Psi$.
\begin{lem}\label{4.4}
Under the same assumption as in Proposition \ref{4.2},
there holds that 
\begin{equation}\label{higher1}
\sup_{t\in[0,T]}\|\nabla^k\Psi(t)\|_{0,\beta}^2
\leq C(\|\Psi_0\|_{k,\beta}^2
+\sup_{t\in[0,T]}\|\Psi(t)\|_{k-1,\beta}^2+|\phi_b|)
\quad \text{$k = 1,\ldots,r$},
\end{equation}
where $C>0$ is a constant depending on $\beta$ but independent of $\phi_b$ and $t$.
\end{lem}
\begin{proof}
To obtain \eqref{higher1}, we estimate each terms on the left hand side of \eqref{hes4} from below separately.
The second term on the left hand side is nonnegative thanks to \eqref{super1}.
The third term is bounded from below as
\begin{align*}
{}&
-\beta\sum_{|\bm{a}|=k}\int_{\Omega} e^{\beta x_1}
\langle A^1[V]\partial_x^{\bm{a}}\Psi,\partial_x^{\bm{a}}\Psi \rangle \,dx
\notag\\
& \geq -\beta\sum_{|\bm{a}|=k}\int_{\Omega} e^{\beta x_1}
\langle A^1[V_{+}]\partial_x^{\bm{a}}\Psi,\partial_x^{\bm{a}}\Psi \rangle \,dx
-C(N_{r,\beta}(T)+|\phi_b|)\|\nabla^k\Psi\|_{0,\beta}^2
\notag\\
& \geq c\beta \|\nabla^{k}\Psi\|_{0,\beta}^2
-C(N_{r,\beta}(T)+|\phi_b|)\|\nabla^k\Psi\|_{0,\beta}^2,
\end{align*}
where we have used \eqref{Bohm1} in deriving the last inequality.

Next we deal with the right hand side of \eqref{hes4}.
First it follows from \eqref{ses1}, $\beta\leq \alpha/2$, and $M\in H^\infty(\Omega)$ that
\begin{equation}\label{h2}
\|\bm{h}\|_{r,\beta} \leq C|\phi_b|. 
\end{equation}
We also see from \eqref{eqiv0}, \eqref{ses1}, and \eqref{ellineq4}
with the aid of the general inequalities \eqref{A2} and \eqref{A4} that
\begin{equation}\label{I1}
\|I^{\bm{a}}\|_{0,\beta} 
\leq C(N_{r,\beta}(T)+|\phi_b|)\|\nabla^k\Psi\|_{0,\beta}
+C\|\Psi\|_{k-1,\beta},
\end{equation}
where we have also used the fact that the component of $(A^0)^{-1} A^j$ is a linear combination of $\Psi$ and smooth functions independent of $\Psi$ in applying \eqref{A2} and \eqref{A4}. 
Now the right hand side of \eqref{hes4} can be estimated as follow:
\begin{align*}
&\sum_{|\bm{a}|=k}\int_{\Omega} e^{\beta x_1} \left(
\langle \{\partial_{t}(A^0[{V}])\}\partial_x^{\bm{a}}\Psi,\partial_x^{\bm{a}}\Psi \rangle
+\sum_{j=1}^3 \langle \{\partial_{x_j}(A^j[{V}])\}\partial_x^{\bm{a}}\Psi,\partial_x^{\bm{a}}\Psi \rangle
\right) \,dx
\notag\\
&\quad +\sum_{|\bm{a}|=k}\int_{\Omega} e^{\beta x_1} 
\left(\partial_x^{\bm{a}}\bm{h}\cdot\partial_x^{\bm{a}} \bm{\eta}
+\langle A^{0}[V] I^{\bm{a}}, \partial_x^{\bm{a}}\Psi \rangle \right) \,dx
\notag\\
&\leq C(N_{r,\beta}(T)+|\phi_b|+\mu)\|\nabla^k\Psi\|_{0,\beta}^2
+C_{\mu}\|\Psi\|_{k-1,\beta}^2+C|\phi_b|,
\end{align*}
where we have also used \eqref{ses1},  \eqref{ellineq4}, \eqref{DtPsi1}, and Schwarz's inequality.

Substituting all the above estimates  into \eqref{hes4} and 
taking $\mu$, $N_{r,\beta}(T)$, and $|\phi_b|$ suitably small, we have
\begin{align*}
 \frac{d}{dt}\sum_{|\bm{a}|=k}\int_\Omega e^{\beta x_1} \langle A^0[V]\partial_x^{\bm{a}}\Psi,\partial_x^{\bm{a}}\Psi \rangle \,dx
+c\|\nabla^k \Psi\|_{0,\beta}^2
\leq C\|\Psi\|_{k-1,\beta}^2+C|\phi_b|.
\end{align*}
Then we multiply this by $e^{\tilde{c} t}$,  integrate over $[0,t]$,
and let $\tilde{c}>0$ be small enough to obtain
\begin{align*}
{}&
e^{\tilde{c} t}\|\nabla^k\Psi(t)\|_{0,\beta}^2
+c\int_0^t e^{\tilde{c}\tau}  \|\nabla^k\Psi(\tau)\|_{0,\beta}^2 \,d\tau
\\
&\leq C\|\Psi_0\|_{k,\beta}^2
+C\int_0^t e^{\tilde{c}\tau}(
\|\Psi(\tau)\|_{k-1,\beta}^2+|\phi_b|) \,d\tau
\\
&\leq C\|\Psi_0\|_{k,\beta}^2
+C\left(\sup_{t \in [0,T]}\|\Psi(t)\|_{k-1,\beta}^2+|\phi_b|\right)
\frac{1}{\tilde{c}}\left(e^{\tilde{c} t}-1\right).
\end{align*}
This immediately gives \eqref{higher1}.
\end{proof}

\subsection{Completion of a priori estimate}\label{S4.4}
We now complete the derivation of the a priori estimate \eqref{apes2}.
\begin{proof}[Proof of Proposition \ref{4.2}]
We begin by proving that
\begin{equation}\label{apes0}
\sup_{t\in[0,T]} \|\Psi(t)\|_{r,\beta}^2 \leq C(\|\Psi_0\|_{r,\beta}^2+|\phi_b|).  
\end{equation}
Substituting \eqref{higher1} with $k=1$ into the right hand side of \eqref{basic1} and 
taking $N_{r,\beta}(T)$ and $|\phi_b|$ sufficiently small, we have
$\sup_{t\in[0,T]}\|\Psi(t)\|_{0,\beta}^2 \leq C(\|\Psi_0\|_{1,\beta}^2+|\phi_b|)$.
Then substituting this into the right hand side of \eqref{higher1} with $k=1$ leads to
$\sup_{t\in[0,T]}\|\Psi(t)\|_{1,\beta}^2 \leq C(\|\Psi_0\|_{1,\beta}^2+|\phi_b|)$.
Furthermore, the induction by using this and \eqref{higher1} yields \eqref{apes0}.

We can complete the proof of \eqref{apes2} by showing that
\begin{gather}
\|\Psi_{t}(t)\|_{l,\beta}^2 \leq C\|\Psi(t)\|_{l+1,\beta}^2+C|\phi_b|^2,
\quad \text{for $l=0,\ldots,r-1$.}
\label{hypineq3}
\end{gather}
Indeed, this with \eqref{ellineq4} and \eqref{apes0} immediately gives \eqref{apes2}.
Let us prove \eqref{hypineq3} for $l=0$. 
Multiply \eqref{re1} by $e^{\beta x_1/2}(A^{0})^{-1}$, 
take the $L^2$-norm, and use \eqref{ellineq4} and \eqref{h2} to obtain
\begin{align*}
\|\Psi_{t}\|_{0,\beta}
&= \left\|
\sum_{j=1}^3 (A^0 )^{-1} A^j \Psi_{x_{j}}
-\begin{bmatrix}
0 \\ m^{-1} \nabla \sigma \\0
\end{bmatrix}
- \tilde{B}\Psi
-\begin{bmatrix}
0 \\ m^{-1}\bm{h}\\0
\end{bmatrix}
\right\|_{0,\beta}
\leq  C(\|\Psi(t)\|_{1,\beta}+|\phi_b|).
\end{align*}
Similarly, we deduce \eqref{hypineq3} for all $l\geq1$ by using \eqref{eqiv0}, \eqref{hes1}, \eqref{h2}, and \eqref{I1}.
The proof is complete.
\end{proof}

\section{Construction of stationary solutions}\label{S5}
This section is devoted to the construction of 
solutions $(\Psi^s,\sigma^s)$ of
the associated stationary problem of \eqref{re0}.
It is to be expected from the bound \eqref{apes1} of time-global solutions $(\Psi,\sigma)$
that these global solutions may converge to 
some functions as $t$ tends to infinity.
Therefore, in subsection \ref{S5.1}, we define an sequence $\{(\Psi^k,\sigma^k)\}_{k=0}^\infty$ by
$(\Psi^k,\sigma^k)(t,x):=(\Psi,\sigma)(t+kT_*,x)$ for any $T^*>0$, and show that
this sequence converges to a time-periodic solution with a period $T^*$ 
to the problem of the equations \eqref{re1} and \eqref{re2} 
with the boundary conditions \eqref{re4} and \eqref{re5}.
By the arbitrary of the period $T^{*}$, 
it can be concluded in subsection \ref{S5.2} that 
the periodic solution is independent of time and 
thus the desired stationary solution.

\subsection{Time-periodic solutions}\label{S5.1}
\subsubsection{Uniqueness}\label{S5.1.1}
We begin by studying the uniqueness of time-periodic solutions 
of the problem \eqref{re1}--\eqref{re5} in the solution space 
\[
{\cal X}^r_{\beta}([0,T^*]):=
\left[L^\infty([0,T^*];H^{r}_{\beta}(\Omega)) 
\cap  W^{1,\infty}([0,T^*];H^{r-1}_{\beta}(\Omega))\right]^5
\times
C([0,T^*];H^{r+1}_{\beta}(\Omega)).
\]
The uniqueness is summarized in the next proposition.
\begin{pro}\label{5.1}
Let $u_+$ and $\theta_{+}$ satisfy \eqref{Bohm1} and \eqref{asp1'}.
For $\beta>0$ being in Theorem \ref{4.1},
there exists $\delta_0>0$ such that 
if a time-periodic solution $(\Psi^*,\sigma^*) \in {\cal X}^3_{\beta}([0,T^*])$
with a period $T^*>0$ of the problem \eqref{re1}--\eqref{re5}
exists and satisfies the following inequality, then it is unique:
\begin{equation}\label{uniasp1}
\sup_{t\in[0,T^*]}(\|\Psi^*(t)\|_{3,\beta}
+\|\Psi^*_{t}(t)\|_{2,\beta}+\|\sigma^*(t)\|_{4,\beta})
+|\phi_b| \leq \delta_0.
\end{equation}
\end{pro}

Let $(\Psi,\sigma)$ and $(\Psi^*,\sigma^*)$ 
be time-periodic solutions with \eqref{uniasp1}, where
$\Psi={}^{t}(\psi,\bm{\eta},\zeta)$ and $\Psi^*={}^{t}(\psi^*,\bm{\eta}^*,\zeta^{*})$.
It is straightforward to see 
that $\bar{\Psi}={}^{t}(\bar{\psi},\bar{\bm{\eta}},\bar{\zeta}):={}^{t}(\psi-\psi^*,\bm{\eta}-\bm{\eta}^*,\zeta-\zeta^{*})$ 
and $\bar{\sigma}:=\sigma-\sigma^*$ satisfy
\begin{subequations}\label{unieq0}
\begin{gather}
\begin{aligned}
&A^0[V] \partial_t\bar{\Psi}
+\sum_{j=1}^3A^j[V]\partial_{x_j}\bar{\Psi}
\\
&=\begin{bmatrix}
0 \\ \nabla \bar{\sigma} \\ 0
\end{bmatrix}
+B[V,\tilde{V}',\nabla M]\bar{\Psi}
-\sum_{j=1}^3 A^0[V] \left\{ ((A^{0})^{-1}A^j)[V] - ((A^{0})^{-1}A^j)[V^{*}] \right\} \partial_{x_j}\Psi^{*},
\end{aligned}
\label{unieq1}
\\
\Delta \bar{\sigma}-\bar{\sigma}=\bar{\psi}
+g_0[\psi,\tilde{v}]-g_0[\psi^*,\tilde{v}]
+g_1[\sigma,\tilde{\phi}]-g_1[\sigma^*,\tilde{\phi}],
\label{unieq2}
\\
 \lim_{|x|\to\infty}(\bar{\Psi},\bar{\sigma})(t,x)=0, 
\label{unieq3}
\\
 \bar{\sigma}(t,M(x'),x')=0,
\label{unieq4}
\end{gather}
\end{subequations}
where $V=\tilde{V}+\Psi$ and $V^{*}=\tilde{V}+\Psi^{*}$.
We note that the system in \eqref{unieq0} is parallel to the system \eqref{re0}.
Indeed, the essential difference between \eqref{re1} and \eqref{unieq1} is only the rightmost.
Furthermore, for the difference between \eqref{re2} and \eqref{unieq2}, 
the terms $g_0[\psi,\tilde{v}]$, $g_1[\phi,\tilde{\phi}]$, and $g_2[\tilde{\phi}',\nabla M]$ 
are only replaced by $g_0[\psi,\tilde{v}]-g_0[\psi^*,\tilde{v}]$, 
$g_1[\sigma,\tilde{\phi}]-g_1[\sigma^*,\tilde{\phi}]$, and zero, respectively.

To show Proposition \ref{5.1}, we use the estimates of $\bar{\sigma}$ in the next lemma.
\begin{lem}\label{uniell1}
Under the same assumption as in Proposition \ref{5.1}, there holds that
\begin{gather}
\|\bar{\sigma}(t)\|_{1,\beta}^2 
\leq (1+D\beta^2+C\delta_0) \|\bar{\psi}(t)\|_{0,\beta}^2,
\label{uniellineq1}\\
\|\bar{\sigma}(t)\|_{2,\beta}^2 \leq C \|\bar{\psi}(t)\|_{0,\beta}^2,
\label{uniellineq4}
\end{gather}
where $D$ is the same positive constant being in Lemma \ref{ell1}, 
and $C$ is a positive constant independent of $\beta$, $\phi_b$, and $t$.
\end{lem}
\begin{proof}
We can show this lemma by using \eqref{unieq2} parallel to \eqref{re2}
in much the same way as in the proof of Lemma \ref{ell1}.
\end{proof}

We are now at a position to show Proposition \ref{5.1}.

\begin{proof}[Proof of Proposition \ref{5.1}]
It suffices to prove $\bar{\Psi}=0$,
since this and \eqref{uniellineq4} lead to $\bar{\sigma}=0$.
To this end, we only need to show
\begin{gather}
\int_0^{T^*} \|\bar{\Psi}(\tau)\|_{0,\beta}^2\,d\tau
\leq C\delta_0 \int_0^{T^*} \|\nabla\bar{\Psi}(\tau)\|_{0,\beta}^2\,d\tau,
\label{unibasic1}\\
\int_0^{T^*} \|\nabla \bar{\Psi}(\tau)\|_{0,\beta}^2\,d\tau
\leq C\int_0^{T^*} \|\bar{\Psi}(\tau)\|_{0,\beta}^2\,d\tau.
\label{unihigher1}
\end{gather}
In fact, one can deduce $\bar{\Psi}=0$ 
by substituting \eqref{unihigher1} into the right hand side of \eqref{unibasic1}
and taking $\delta_0$ sufficiently small.

From now on we show \eqref{unibasic1}.
We recall that \eqref{unieq1} is parallel to \eqref{re1}, and then
derive the following equality parallel to \eqref{bes7}:
\begin{align}
{}&
\frac{d}{dt} \int_\Omega e^{\beta x_1}
\left(\bar\psi^2+
\langle  A^0[{V}]\bar\Psi,\bar\Psi\rangle+m(\nabla \cdot  \bar{\bm{\eta}})^2+R\theta |\nabla \bar\psi|^2+\frac{R}{(\gamma-1)\theta}|\nabla\bar\zeta|^2\right)\,dx
\notag\\
&
\quad+\sum_{j=1}^3\int_{\partial\Omega} e^{\beta M(x')}
\langle n_j A^j[V]\bar\Psi,\bar\Psi \rangle \,ds
+\sum_{j=1}^3\int_{\partial\Omega} e^{\beta M(x')}n_ju_{j}\bar\psi^2\,ds
\notag\\
&
\quad+\int_{\partial\Omega} e^{\beta M(x')}\left\langle\!
F[V,\bm{n}]
\begin{bmatrix}
\nabla \cdot  \bar{\bm{\eta}}
\\
\nabla \bar\psi
\\
\nabla \bar\zeta
\end{bmatrix},
\begin{bmatrix}
\nabla \cdot  \bar{\bm{\eta}}
\\
\nabla \bar\psi
\\
\nabla \bar\zeta
\end{bmatrix}
\!\right\rangle\,ds
\notag\\
&
\quad-\beta\int_{\Omega} e^{\beta x_1}
(\langle A^1[V]\bar\Psi,\bar\Psi \rangle
+u_{1}\bar\psi^2
-2\bar\sigma\bar\eta_1)\,dx
-\beta\int_{\Omega} e^{\beta x_1}
\left\langle\!
F^1[{V}]
\begin{bmatrix}
\nabla \cdot  \bar{\bm{\eta}}
\\
\nabla \bar\psi
\\
\nabla \bar\zeta
\end{bmatrix},
\begin{bmatrix}
\nabla \cdot  \bar{\bm{\eta}}
\\
\nabla \bar\psi
\\
\nabla \bar\zeta
\end{bmatrix}
\!\right\rangle\,dx
\notag\\
&=\int_{\Omega} e^{\beta x_1} \sum_{i=1}^{5} \bar{\mathcal R}_{i}\,dx,
\label{unibes10}
\end{align}
where the $7 \times 7$ symmetric matrices $F$ and $F^{1}$ are defined in \eqref{DefF} and \eqref{DefF1}, respectively, \footnote{We note that $\bar{\mathcal{R}}_1,\ldots,\bar{\mathcal{R}}_5$ correspond to ${\mathcal{R}}_1,\ldots,{\mathcal{R}}_5$ in subsection \ref{S4.2}, respectively.}and $\bar{\mathcal{R}}_i$ for $i=1,\ldots,5$ is defined as
\begin{align*}
\bar{\mathcal{R}}_1
&:=\langle \{(A^0[{V}])_{t}\}\bar{\Psi},\bar{\Psi} \rangle+\sum_{j=1}^3 \langle \{(A^j[{V}])_{x_{j}}\}\bar{\Psi},\bar{\Psi} \rangle +2\langle B\bar{\Psi}, \bar{\Psi}\rangle
\\
& \quad -\langle\sum_{j=1}^3 A^0[V] \left\{ ((A^{0})^{-1}A^j)[V] - ((A^{0})^{-1}A^j)[V^{*}] \right\} \Psi^{*}_{x_{j}}, \bar{\Psi}\rangle,
\\
\bar{{\mathcal R}}_2
&:=m\sum_{j=1}^3 (u_{j})_{x_j}(\nabla \cdot  \bar{\bm{\eta}})^2
 -2\Big(m\nabla u\cdot\nabla  \bar{\bm{\eta}}+m\nabla  \bar{\bm{\eta}}\cdot\nabla {\bm{\eta}}^*+m \bar{\bm{\eta}}\cdot
\nabla(\nabla {\bm{\eta}}^*)\Big)(\nabla \cdot   \bar{\bm{\eta}})
\\
&
 \quad -2\left( R\sum_{i=1}^3 \theta_{x_i} \bar\psi_{x_i}+R\sum_{i,j=1}^3  \bar{{\zeta}}_{x_i}\psi^*_{x_i}+R\bar\zeta\Delta\psi^*
-\sum_{i=1}^3\{(B\Psi)_{1+I}\}_{x_i}\right)(\nabla \cdot   \bar{\bm{\eta}}),
\\
\bar{{\mathcal R}}_3
&:=R\theta_t|\nabla \bar\psi|^2+
\sum_{j=1}^3 \left\{R \theta u_{j} \right\}_{x_j}|\nabla \bar\psi|^2 +2R(\nabla \cdot  \bar{\bm{\eta}}) (\nabla \theta \cdot \nabla \bar\psi)
\\
&  \quad -2\Big( R \theta \nabla \bm{u} \cdot\nabla \bar\psi +R \theta \nabla \bar{\bm{\eta}} \cdot\nabla \psi^*
+R \theta  \bar{\bm{\eta}} \cdot\nabla \nabla\psi^*-R\theta \nabla(\tilde{B}\Psi)_1\Big)\cdot \nabla \bar\psi,
\\
\bar{{\mathcal R}}_4
&:=\left\{\frac{R}{(\gamma-1)\theta }\right\}_t|\nabla \bar\zeta|^2+
\sum_{j=1}^3 \left\{\frac{Ru_{j}}{(\gamma-1)\theta}\right\}_{x_j}|\nabla \bar\zeta|^2
-2\frac{R}{\theta}\left(\frac{\nabla \bm{u}\cdot \nabla \bar\zeta}{\gamma-1}
+\frac{\nabla \bar{\bm{\eta}}\cdot \nabla \zeta^*}{\gamma-1} \right) \cdot \nabla \bar\zeta
\\
& \quad -2\frac{R}{\theta}\left( \frac{\bar{\bm{\eta}}\cdot \nabla\nabla \zeta^*}{\gamma-1}
+(\nabla\cdot \bar{\bm{\eta}})\nabla\theta
+(\nabla\cdot {\bm{\eta}}^*)\nabla\bar\zeta
+(\nabla\nabla\cdot {\bm{\eta}}^*)\bar\zeta
-\frac{\nabla(\tilde{B}\Psi)_5}{\gamma-1}\right)\cdot \nabla \bar\zeta,
\\
\bar{{\mathcal R}}_5
&:=2(\nabla\cdot \bm{u}) \bar{\psi}^{2}
+ 2\bar{\psi}(B\bar{\Psi})_1
-2\bar{\psi}(\bar{\bm{\eta}}\cdot\nabla \psi^*) 
\\
&\quad +2(g_0[\psi,\tilde{v}]-g_0[\psi^*,\tilde{v}]
+g_1[\sigma,\tilde{\phi}]-g_1[\sigma^*,\tilde{\phi}])(\nabla\cdot\bar{\bm{\eta}}).
\end{align*}

Then we notice that the left hand side of \eqref{unibes10} has much the same form as that of \eqref{bes7}.
Therefore, the second and third terms are nonnegative and so negligible 
if $\delta_0$ is sufficiently small.
Furthermore, with the aid of \eqref{uniellineq1}, 
the fourth and fifth terms are bounded from below as
\begin{align*}
{}&
-\beta\int_{\Omega} e^{\beta x_1}
(\langle A^1[{V}]\bar{\Psi},\bar{\Psi} \rangle
+(\eta_1+\tilde{u})\bar{\psi}^2
-2\bar{\sigma}\bar{\eta}_1)\,dx
\geq (\beta d - 2\sqrt{D}\beta^2+\mu+C_{\mu} \delta_0)\|\bar{\Psi}\|_{0,\beta}^2,
\\
{}&
-\beta\int_{\Omega} e^{\beta x_1}
\left\langle\!
F^1[V]
\begin{bmatrix}
\nabla \cdot  \bm{\eta}
\\
\nabla \psi
\\
\nabla \zeta
\end{bmatrix},
\begin{bmatrix}
\nabla \cdot  \bm{\eta}
\\
\nabla \psi
\\
\nabla \zeta
\end{bmatrix}
\!\right\rangle\,dx
\geq \beta d\|(\nabla \cdot \bm{\eta},\nabla \psi,\nabla\zeta)\|_{0,\beta}^2
-C\delta_{0}\|\Psi\|_{1,\beta}^2,
\end{align*}
where $d$ and $D$ are the same positive constants as in \eqref{defbeta1}.
Let us also estimate the right hand side of \eqref{unibes10}, but we must be careful to handle the terms 
having $\psi^*$ and $\bm{\eta}^*$ in $\bar{\mathcal R}_{i}$,
since some of these include the second-order derivatives.
Using \eqref{ses1} and \eqref{uniasp1}, we estimate $\bar{\mathcal R}$ as
\begin{equation*}
|\bar{\mathcal R}_{i}|\leq
C\delta_0|(\bar{\Psi},\bar{\sigma},\nabla\bar{\Psi})|^2
+C|\bar{\Psi}||\nabla\bar{\Psi}||\nabla^2{\Psi^*}| \quad \text{for $i=1,\ldots,5$}.
\end{equation*}
Then Sobolev's and H\"older's inequalities give
\begin{align}\label{unibes9}
\left|\int_{\Omega} e^{\beta x_1}\bar{\mathcal R}_{i} \,dx\right|
\leq C\delta_0(\|\bar{\Psi}\|_{1,\beta}^2+\|\bar{\sigma}\|_{1,\beta}^2)
\leq C\delta_0\|\bar{\Psi}\|_{1,\beta}^2,
\end{align}
where we have also used \eqref{uniellineq4} in deriving the last inequality.

Now we substitute these inequalities into \eqref{unibes10}, use \eqref{defbeta1}, 
and take $\mu$ and $\delta_0$ small enough to obtain
\begin{multline}
\frac{d}{dt}  \int_\Omega e^{\beta x_1}
\left(\bar\psi^2+
\langle  A^0[{V}]\bar\Psi,\bar\Psi\rangle+m(\nabla \cdot  \bar{\bm{\eta}})^2+R\theta |\nabla \bar\psi|^2+\frac{R}{(\gamma-1)\theta}|\nabla\bar\zeta|^2\right)\,dx
\\
+c\beta\|(\bar{\Psi},\nabla \bar{\psi},\nabla \cdot \bar{\bm{\eta}},\nabla\bar\zeta)\|_{0,\beta}^2
\leq C\delta_0\|\nabla\bar{\Psi}\|_{0,\beta}^2.
\label{unibes12}
\end{multline}
Then integrating this over $[0,T^*]$ and using the periodicity of solutions, 
we conclude \eqref{unibasic1}.

Let us complete the proof by showing \eqref{unihigher1}.
Multiply \eqref{unieq1} by $(A^{0}[V])^{-1}$, 
apply $\partial_x^{\bm{a}}$ with $|\bm{a}|=1$ to the result,
multiply it by $A^{0}[V]$, 
take an inner product of the result with $2e^{\beta x_1}\partial_x^{\bm{a}}\bar{\Psi}$,
and sum up the results for $\bm{a}$ with $|\bm{a}|=1$.
Then integrating the resultant equality over $\Omega$ 
and applying Gauss's divergence theorem, we have
\begin{multline}
\frac{d}{dt} \sum_{|\bm{a}|=1}\int_\Omega e^{\beta x_1} 
\langle A^0[{V}]\partial_x^{\bm{a}}\bar{\Psi},\partial_x^{\bm{a}}\bar{\Psi} \rangle \,dx
+\sum_{|\bm{a}|=1}\sum_{j=1}^3\int_{\partial\Omega} e^{\beta M(x')}
\langle n_j A^j[{V}]\partial_x^{\bm{a}}\bar{\Psi},\partial_x^{\bm{a}}\bar{\Psi} \rangle\,ds
\\
-\beta\sum_{|\bm{a}|=1}\int_{\Omega} e^{\beta x_1}
\langle A^1[{V}]\partial_x^{\bm{a}}\bar{\Psi},\partial_x^{\bm{a}}\bar{\Psi} \rangle\,dx
= \sum_{|\bm{a}|=1}\int_{\Omega} e^{\beta x_1}  \bar{\mathcal{R}}_\alpha  \,dx,
\label{unihes4}
\end{multline}
where
\begin{align*}
\bar{\mathcal{R}}_\alpha &:=  \langle \{\partial_{t}(A^0[{V}])\}\partial_x^{\bm{a}}\bar{\Psi},\partial_x^{\bm{a}}\bar{\Psi} \rangle+
\sum_{j=1}^3 \langle \{\partial_{x_j}(A^j[{V}])\}\partial_x^{\bm{a}}\bar{\Psi},\partial_x^{\bm{a}}\bar{\Psi} \rangle 
\notag\\
&\quad +2(\nabla\partial_x^{\bm{a}}\bar{\sigma})\cdot\partial_x^{\bm{a}}\bar{\bm{\eta}} 
+2\langle A^{0}[V]\partial_x^{\bm{a}}(\tilde{B}\bar{\Psi}), \partial_x^{\bm{a}}\bar{\Psi} \rangle 
\notag\\
&\quad -2\sum_{j=1}^3\langle A^{0}[V] \{\partial_x^{\bm{a}}((A^{0}[V])^{-1}A^j[{V}])\}\partial_{x_j}\bar{\Psi},\partial_x^{\bm{a}}\bar{\Psi} \rangle 
\notag\\
&\quad -2\langle A^0[V]  \partial_{x}^{\bm{a}} \sum_{j=1}^3   \left\{ ((A^{0})^{-1}A^j)[V] - ((A^{0})^{-1}A^j)[V^{*}] \right\} \partial_{x_j}\Psi^{*},\partial_x^{\bm{a}}\bar{\Psi} \rangle.
\end{align*}

We notice that the left hand side of the equality \eqref{unihes4}
has the same form as that of \eqref{hes4}.
Therefore, the second term on the left hand side is nonnegative.
The third term is bounded from below as
\begin{align*}
-\beta\sum_{|\bm{a}|=1}\int_{\Omega} e^{\beta x_1}
\langle A^1[{V}]\partial_x^{\bm{a}}\bar{\Psi},\partial_x^{\bm{a}}\bar{\Psi} \rangle \,dx
\geq c\beta \|\nabla\bar{\Psi}\|_{0,\beta}^2
-C\delta_0\|\nabla\bar{\Psi}\|_{0,\beta}^2.
\end{align*}
On the other hand, the right hand side of \eqref{unihes4} can be estimated similarly to \eqref{unibes9} as follows:
\begin{gather*}
\sum_{|\bm{a}|=1}\int_{\Omega} e^{\beta x_1}  \bar{\mathcal{R}}_\alpha  \,dx
\leq (C\delta_0+\mu)\|\bar{\Psi}\|_{1,\beta}^2
+C_{\mu}\|\nabla^2\bar{\sigma}\|_{0,\beta}^2
\leq (C\delta_0+\mu)\|\nabla\bar{\Psi}\|_{0,\beta}^2
+C_{\mu}\|\bar{\Psi}\|_{0,\beta}^2.
\end{gather*}
We substitute these estimates into \eqref{unihes4} and
let $\mu$ and $\delta_0$ be sufficiently small to obtain
\begin{equation}\label{unihes5}
\frac{d}{dt} \sum_{|\bm{a}|=1}\int_\Omega e^{\beta x_1} 
\langle A^0[{V}]\partial_x^{\bm{a}}\bar{\Psi},\partial_x^{\bm{a}}\bar{\Psi} \rangle \,dx
+c\|\nabla\bar{\Psi}\|_{0,\beta}^2
\leq C\|\bar{\Psi}\|_{0,\beta}^2. 
\end{equation}
Then integrating this over $[0,T_*]$ and 
using the periodicity of solutions, we conclude \eqref{unihigher1}.
\end{proof}

\subsubsection{Existence}\label{S5.1.2}
In this subsection, we establish the existence of time-periodic solutions of the problem \eqref{re1}--\eqref{re5}.
Specifically, we show the following proposition:
\begin{pro}\label{5.3}
Let $r \geq 3$ as well as $u_+$ and $\theta_+$ satisfy \eqref{Bohm1} and \eqref{asp1'}.
For $\beta>0$ being in Theorem \ref{4.1} and any $T^*>0$,
there exists a constant $\delta>0$ independent of $T^*$ such that 
if $|\phi_b| \leq \delta$, then the problem \eqref{re1}--\eqref{re5} 
has a time-periodic solution $(\Psi^*,\sigma^*)\in {\mathcal X}^r_\beta([0,T^*])$
with a period $T^*>0$. Furthermore, it satisfies
\begin{equation*}
\sup_{t\in[0,T^*]}(\|\Psi^*(t)\|_{r,\beta}
+\|\partial_t \Psi^*(t)\|_{r-1,\beta}+\|\sigma^*(t)\|_{r+1,\beta})
\leq C|\phi_b|^{1/2},
\end{equation*}
where $C>0$ is a constant independent of $T^*$.
\end{pro}

For the construction of time-periodic solutions, we define
\[
(\Psi^k,\sigma^k)(t,x):=(\Psi,\sigma)(t+kT^*,x)
\quad \text{for $k=1,2,3,\ldots$,} 
\]
where $(\Psi,\sigma)$ is the time-global solution in Theorem \ref{4.1}
and $\Psi^k$ denotes ${}^t(\psi^k,\bm{\eta}^k,\zeta^k)$.
To show Proposition \ref{5.3}, it suffice to show Lemma \ref{5.2} below.
Indeed, Proposition \ref{5.3} follows from much the same argument in the proof of \cite[Proposition 5.4]{MM1} with the aid of Lemma \ref{5.2}. 
In the process, we can know that 
\begin{equation}\label{converge1}
(\Psi^k,\sigma^k) \to (\Psi^*,\sigma^*) \quad \text{in} \ \
\left[\bigcap_{i=0}^1C^i([0,T^*];H^{r-i-1}_\beta(\Omega))\right]^5
\times C([0,T^*];H^{r+1}_\beta(\Omega)),
\end{equation}
where $(\Psi^{*},\sigma^{*})$ is the time-periodic solution in Proposition \ref{5.3}.

\begin{lem}\label{5.2}
Let $u_+$ and $\theta_+$ satisfy \eqref{Bohm1} and \eqref{asp1'}.
For $\beta>0$ being in Theorem \ref{4.1} and any $T^*>0$,
there exist $\lambda>0$ and $C>0$ independent of $k$ and $T^*$ such that
\begin{equation}\label{exiapes0}
\|(\Psi-\Psi^k)(t)\|_{1,\beta}
+\|(\sigma-\sigma^k)(t)\|_{2,\beta} \leq Ce^{-\lambda t}
\quad \text{for $k=1,2,3,\ldots$.}
\end{equation}
\end{lem}
\begin{proof}
We recall that the time-global solution in Theorem \ref{4.1} satisfies \eqref{apes1}.
Therefore, by the same method as in the derivations of \eqref{unibes12} and \eqref{unihes5}, it is seen that
\begin{align*}
{}&
\frac{d}{dt} \int_\Omega e^{\beta x_1}
\left(|\psi-\psi^k|^2+
\langle  A^0[{V}](\Psi-\Psi^k),(\Psi-\Psi^k)\rangle \right)\,dx
\\
&\quad +\frac{d}{dt} \int_\Omega e^{\beta x_1}
\left(m|\nabla \cdot   (\bm{\eta}-\bm{\eta}^k)|^2 + R\theta |\nabla (\psi-\psi^k)|^2+\frac{R}{(\gamma-1)\theta}|\nabla(\zeta-\zeta^k)|^2\right)\,dx
\\
&\quad +c\beta\|(\Psi-\Psi^k,\nabla \cdot (\bm{\eta}-\bm{\eta}^k),\nabla(\psi-\psi^k),\nabla(\zeta-\zeta^k))\|_{0,\beta}^2\\
&\leq C(\|\Psi_0\|_{r,\beta}+|\phi_b|^{1/2})\|\nabla(\Psi-\Psi^k)\|_{0,\beta}^2
\end{align*}
and
\begin{gather*}
\frac{d}{dt} \sum_{|\bm{a}|=1}\int_\Omega e^{\beta x_1} 
\langle A^0[{V}]\partial_x^{\bm{a}}(\Psi-\Psi^k),\partial_x^{\bm{a}}(\Psi-\Psi^k) \rangle \,dx
+c\beta\|\nabla(\Psi-\Psi^k)\|_{0,\beta}^2
\leq C\|\Psi-\Psi^k\|_{0,\beta}^2. 
\label{exihes5}
\end{gather*}
Then we multiply these two by $e^{\tilde{c}t}$, 
integrate the results over $[0,T^*]$, 
and take $\tilde{c}>0$ suitably small to obtain
\begin{align*}
{}& e^{\tilde{c}t}\|(\Psi-\Psi^k)(t)\|_{0,\beta}^2
+\int_0^t e^{\tilde{c}\tau}\|(\Psi-\Psi^k)(\tau)\|_{0,\beta}^2\,d\tau
\\
&\leq C\|(\Psi-\Psi^k)(0)\|_{1,\beta}^2
+C(\|\Psi_0\|_{r,\beta}+|\phi_b|^{1/2})
\int_0^t e^{\tilde{c}\tau}\|\nabla(\Psi-\Psi^k)(\tau)\|_{0,\beta}^2\,d\tau
\end{align*}
and
\begin{align*}
{}& e^{\tilde{c}t}\|\nabla(\Psi-\Psi^k)(t)\|_{0,\beta}^2
+\int_0^t e^{\tilde{c}\tau}\|\nabla(\Psi-\Psi^k)(\tau)\|_{0,\beta}^2\,d\tau
\\
&\leq C\|(\Psi-\Psi^k)(0)\|_{1,\beta}^2
+C\int_0^t e^{\tilde{c}\tau}\|(\Psi-\Psi^k)(\tau)\|_{0,\beta}^2\,d\tau.
\end{align*}
From these two and \eqref{apes1},
we have the estimate of $\Psi-\Psi^k$ in \eqref{exiapes0}
by taking $\|\Psi_0\|_{r,\beta}$ and $|\phi_b|$ suitably small again if necessary.
Now it remains to obtain the estimate of $\sigma-\sigma^k$ in \eqref{exiapes0}.
The same proof as Lemma \ref{uniell1} works for $\sigma-\sigma^k$ and thus 
$\|\sigma-\sigma^k\|_{2,\beta}\leq C\|\psi-\psi^k\|_{0,\beta}$ holds.
This immediately completes the proof.
\end{proof}

\subsection{Stationary solutions}\label{S5.2}
We complete the proof of Theorem \ref{2.1} stating the unique existence of stationary solutions.

\begin{proof}[Proof of Theorem \ref{2.1}]
By much the same argument as in the proof of \cite[Theorem 3.1]{MM1},
we see from the arbitrary of period and the uniqueness of the time-periodic solution that the time-periodic solution $(\Psi^*,\sigma^*)$ in Proposition \ref{5.3} is time-independent.
Hence $(\Psi^s,\sigma^s)=(\Psi^*,\sigma^*)$ is the desired stationary solution. 
The uniqueness follows from Proposition \ref{5.1}. 
The proof is complete.
\end{proof}
%

\subsection{Stability in the exponential weighted Sobolev space}\label{S5.3}
This subsection is devoted to the completion of the proof of Theorem \ref{2.2}.
Since the time-global solutions of the problem \eqref{re0} has been constructed in Theorem \ref{4.1},
it suffices to show the asymptotic stability of stationary solutions.

\begin{proof}[Proof of Theorem \ref{2.2}]
From Theorem \ref{4.1} and Lemma \ref{5.2}, it is seen that the initial--boundary value problem \eqref{re0} has 
a unique time-global solution satisfying \eqref{apes1} and \eqref{exiapes0}
if $\|\Psi_0\|_{r,\beta}$ and $|\phi_b|$ are small enough.
Passing to the limit $k\to \infty$ in \eqref{exiapes0},
we have $\|(\Psi-\Psi^s,\sigma-\sigma^s)(t)\|_{0,\beta} \leq C e^{-\lambda t}$
thanks to \eqref{converge1} and $(\Psi^s,\sigma^s)=(\Psi^*,\sigma^*)$.
Then this inequality and \eqref{apes1} with the aid of the Gagliardo-Nirenberg inequalities 
give the decay estimate \eqref{decay1}.
The proof is complete.
\end{proof}

\medskip

\noindent
{\bf Acknowledge.} 
This work was supported by JSPS KAKENHI Grant Numbers 21K03308.


\begin{thebibliography}{99}

\bibitem{A.A.} 
\newblock {A. Ambroso},
\newblock \emph{Stability for solutions of a stationary Euler--Poisson problem},
\newblock {Math. Models Methods Appl. Sci.} \textbf{16} (2006), 1817--1837.

\bibitem{AMR}
\newblock {A. Ambroso, F. M\'ehats and P.-A. Raviart},
\newblock \emph{On singular perturbation problems for the nonlinear Poisson equation},
\newblock {Asympt. Anal.} \textbf{25} (2001), 39--91.

\bibitem{Ba}
\newblock {M. Badsi},
\newblock \emph{Linear electron stability for a bi-kinetic sheath model},
\newblock {J. Math. Anal. Appl.} \textbf{453} (2017), 954--972.

\bibitem{BCD1}
\newblock {M. Badsi, M. Campos Pinto, and B. Despr\'es},
\newblock \emph{A minimization formulation of a bi-kinetic sheath},
\newblock {Kinet. Relat. Models} \textbf{9} (2016), 621--656.

\bibitem{D.B.1}
\newblock {D. Bohm},
\newblock \emph{Minimum ionic kinetic energy for a stable sheath},
\newblock {in The characteristics of electrical discharges in magnetic fields},
\newblock 
A. Guthrie and R.K.Wakerling eds., McGraw-Hill, New York, (1949), 77--86.

\bibitem{BT1}
\newblock {R. L. F. Boyd and J. B. Thompson},
\newblock \emph{The Operation of Langmuir Probes in Electro-Negative Plasmas},
\newblock {Proc. R. Soc. Lond. A} \textbf{252} (1959), 102--119.

\bibitem{F.C.1}
\newblock {F. F. Chen},
\newblock Introduction to Plasma Physics and Controlled Fusion,
\newblock 2nd edition, Springer, 1984.



\bibitem{DYZ1}
\newblock {R. Duan, H. Yin, and C. Zhu},
\newblock \emph{A half-space problem on the full Euler-Poisson system}
\newblock {SIAM J. Math. Anal.} \textbf{53} (2021), 6094--6121.


\bibitem{FHS1}
\newblock  {M. Feldman, S.-Y. Ha, and M. Slemrod}, 
\newblock  \emph{A geometric level-set formulation of a plasma-sheath interface},
\newblock  {Arch. Ration. Mech. Anal.} \textbf{178} (2005), 81--123.



\bibitem{GHR1}
\newblock {D. G\'erard-Varet, D. Han-Kwan and F. Rousset}, 
\newblock \emph{Quasineutral limit of the Euler--Poisson system for ions
	in a domain with boundaries},
\newblock {Indiana Univ. Math. J.} \textbf{62} (2013), 359--402.

\bibitem{GHR2}
\newblock {D. G\'erard-Varet, D. Han-Kwan and F. Rousset}, 
\newblock \emph{Quasineutral limit of the Euler--Poisson system for ions
	in a domain with boundaries II},
\newblock {J. \'Ec. polytech. Math.} \textbf{1} (2014), 343--386.

\bibitem{GGPS}
\newblock {E. Grenier, Y. Guo, B. Pausader and M. Suzuki}, 
\newblock \emph{Derivation of the ion equation}, 
\newblock {Quarterly of Applied Mathematics}, \textbf{78} (2020), 305--332.




\bibitem{JKS1}
\newblock {C.-Y. Jung, B. Kwon and M. Suzuki},
\newblock \emph{Quasi-neutral limit for the Euler--Poisson system in the presence of plasma sheaths with spherical symmetry},  
\newblock {Math. Models Methods Appl. Sci.} \textbf{26} (2016), 2369--2392.

\bibitem{JKS20}
C.-Y. Jung, B. Kwon, and M. Suzuki.
\newblock \emph{Quasi-neutral limit for Euler--Poisson system in an annular domain},
\newblock J. Differential Equations \textbf{269} (2020), 8007--8054.

\bibitem{JKS21}
C.-Y. Jung, B. Kwon, and M. Suzuki.
\newblock \emph{On approximate solutions to the Euler--Poisson system with boundary layers},
\newblock Commun. Nonlinear Sci. Numer. Simul. \textbf{96} (2021), Paper No. 105717, 27 pp.

\bibitem{I.L.1} 
\newblock {I. Langmuir}, 
\newblock \emph{The interaction of electron and positive ion space charges in cathode sheaths},
\newblock {Phys. Rev.} \textbf{33} (1929), pp. 954--989.

\bibitem{LL}
\newblock {M. A. Lieberman and A. J. Lichtenberg},
\newblock Principles of Plasma Discharges and Materials Processing, 
\newblock 2$^{nd}$ edition, Wiley-Interscience, 2005.


\bibitem{LS1}
\newblock {H. Liu and M. Slemrod}, 
\newblock \emph{KdV dynamics in the plasma-sheath transition},
\newblock {Appl. Math. Lett.} \textbf{17} (2004), 401--410.


\bibitem{NOS} 
\newblock {S. Nishibata, M. Ohnawa and M. Suzuki},
\newblock \emph{Asymptotic stability of boundary layers 
to the Euler--Poisson equations arising in plasma physics},
\newblock {SIAM J. Math. Anal.} \textbf{44} (2012), 761--790. 


\bibitem{K.R.1}
\newblock {K.-U. Riemann},
\newblock \emph{The Bohm criterion and sheath formation. Initial value problems}, 
\newblock {J. Phys. D: Appl. Phys.} \textbf{24} (1991), 493--518.

\bibitem{K.R.2}
\newblock {K.-U. Riemann},
\newblock \emph{The Bohm Criterion and Boundary Conditions for a Multicomponent System}, 
\newblock {IEEE Trans. Plasma Sci.} \textbf{23} (1995), 709--716.

\bibitem{RD1}
\newblock {K.-U. Riemann and T. Daube}, 
\newblock \emph{Analytical model of the relaxation of a collisionless ion matrix sheath},
\newblock {J. Appl. Phys.} \textbf{86} (1999), 1201--1207. 

\bibitem{SS1} 
\newblock {M. Slemrod and N. Sternberg}, 
\newblock \emph{Quasi-neutral limit for Euler-Poisson system},
\newblock {J. Nonlinear Sci.} \textbf{11} (2001), 193--209. 

\bibitem{M.S.1} 
\newblock {M. Suzuki},
\newblock \emph{Asymptotic stability of stationary solutions 
to the Euler--Poisson equations arising in plasma physics},
\newblock {Kinet. Relat. Models} \textbf{4} (2011), 569--588. 

\bibitem{M.S.2} 
\newblock {M. Suzuki},
\newblock \emph{Asymptotic stability of a boundary layer 
to the Euler--Poisson equations for a multicomponent plasma},
\newblock {Kinet. Relat. Models} \textbf{9} (2016), 587-603.

\bibitem{MM1} 
{M. Suzuki and M. Takayama},
\emph{Stability and existence of stationary solutions to 
the Euler--Poisson equations in a domain with a curved boundary},
\newblock {Arch. Ration. Mech. Anal.} \textbf{239} (2021), 357--387. 

\bibitem{MM2} M. Suzuki and M. Takayama, \emph{The Kinetic and Hydrodynamic Bohm Criterions for Plasma Sheath Formation}, Arch. Ration. Mech. Anal. \textbf{247} (2023), 86.

\bibitem{STZ1} 
\newblock {M. Suzuki, M. Takayama, and K. Z. Zhang},
\newblock \emph{Nonlinear Stability and Instability of Plasma Boundary Layers},
\newblock arXiv:2208.07326.

\bibitem{Va83}
\newblock {A.~Valli}, 
\newblock \emph{Periodic and stationary solutions for compressible Navier-Stokes equations via a stability method}, 
\newblock Ann. Scuola Norm. Sup. Pisa Cl. Sci. (4) \textbf{10} (1983),  607--647. 

\end{thebibliography}
\end{document}